\pdfoutput=1
\documentclass[11pt]{article}
\usepackage{amssymb}
\usepackage{graphicx}

\parskip \baselineskip

\newcommand{\s}{\subseteq}
\newcommand{\es}{\emptyset}
\newcommand{\M}{\mathbb{M}}
\newcommand{\La}{\mathbb{L}}

\newcommand{\Ci}{\mathbb{C}}
\newcommand{\N}{\mathbb{N}}
\newcommand{\Pow}{\mathbb{P}}

\newcommand{\B}{{\cal B}}

\newcommand{\Ra}{\Rightarrow}
\newcommand{\LRa}{\Leftrightarrow}

\newcommand{\LL}{\Lambda}

\begin{document}
\centerline{\bf Modular lattices of finite length (Part A)}
\medskip
\centerline{Marcel Wild}
\bigskip

\centerline{\bf 1. Introduction}

These notes constitute the first  revision of my unpublished 1992 notes [W2] with the same title, a  project originally suggested by the late Gian-Carlo Rota.   Laying dormant for  decades  my love for modular lattices got rekindled only recently [W8]. In the process I realized that  nobody yet has heeded Rota's wish, and so I decided to resume this endeavour.
Among others, my teacher (during [W1] and [HW]) Christian Herrmann would have been a more knowledgeable author. On the other hand, he would have been hard pressed to stay within the (challenging enough) realm of finite length modular lattices.

 {\bf Convention:} If not stated otherwise, all lattices are henceforth silently assumed to be of finite length.
 
For the time being these notes will stay on the ResearchGate and in the arXiv, being updated from time to time.
Whether and where they will eventually be published, remains to be seen. Comments and contributions (potentially leading to co-authorship) are welcome. 
The Section titles are as follows.

\begin{enumerate}
\item Introduction
\item Basic facts
\item Lattice congruences via  coverings
\item Geometric lattices and matroids
\item Partial linear spaces and matroids
\item Enumeration of modular lattices\\
\item Bases of lines: The essentials
\item Calculating the submodule lattice of a finite $R$-module
\item Bases of lines: Localization and two types of cycles\\
\item Existence of R-linear representations, partition representations, and 2-distributivity
\item Classification of the k-linear representations of acyclic modular lattices
\item Axiomatization\\
\item Constructions: Gluings and subdirect products
\item Relatively free structures: From semilattices to distributive to modular lattices

\end{enumerate}

As indicated by the spacing, the planned topics are partitioned in Part A up to Part D,  the most finalised of which being Part A. These parts will appear (and will be updated)  individually on ResearchGate and arXiv. Currently only Part A and B are uploaded. Even within individual Parts the levels of coherency may vary considerably. In particular, several '??' remain to be filled. As to the references, currently the author dominates disproportionally; that will change with time (any hints are welcome).

In each Section the numbering of displayed formulas (or technical statements) restarts anew as (1),(2),.. and so forth. Likewise the numbering of Theorems and Lemmas (thus Lemma 3.1, Theorem 3.2, Theorem 3.3 in Section 3), the numbering of Figures, and the numbering of exercises. Speaking of the latter, apart from probing the reader's understanding, they also serve to unclutter the main text. Solutions to most exercises are planned. At present they are given just below the exercises. Missing solutions are due to their triviality or to the author's procrastination.

{\bf 1.1} Here comes an (incomplete) Section break-up.
Section 2 introduces basic lattice concepts, such as easy consequences of distributivity and semimodularity, closure systems and implications. Section 3 propagates to view lattice congruences mainly in terms of the behaviour of coverings (=2-element intervals). This works particularly well for modular lattices and, so it seems, features novel material. Section 4 glances at matroids and the coupled geometric (:= semimodular and atomistic) lattices $L$. A prominent special case are {\it modular} geometric lattices.
As is well known, they are cryptomorphic to projective spaces. Taking another line of attack (Section 5), projective spaces generalize to partial linear spaces (PLSes). The standard text [BB] says nothing about so-called 'point-splittings' that are  crucial  to make a PLS acyclic, and hence more manageable.
 The diagrams of all modular lattices up to cardinality 12 in Section 6 were kindly provided by Kohonen, see also [K].

Section 7  keeps 'modular' but drops 'atomistic' while attempting to maintain the framework of projective geometry as far as possible. It is based on ground-breaking work of Herrmann, Faigle, Jonsson-Nation, Benson-Conway. The crucial concept generalizing a projective geometry is a certain PLS which was coined a 'base of lines' (BOL) in [HW]. Let $L$ be a modular lattice. The point set of a BOL of $L$ is the set $J(L)$ of join irreducibles, and the lines of the BOL are chosen among those sets $\ell\s J(L)$ that are maximal with respect to  any distinct $p,q\in\ell$ yielding the same join $p+q$. Applying bases of lines it is shown in Section 8 how the submodule lattice $\La(M)$ of a finite $R$-module $M$ can be calculated in a compressed format. Modulo an open question that method extends to the calculation of the sublattice $\La'\s\La(M)$ generated by arbitrary members of $\La(M)$. It remains to compare this with the wholly different approach in [LMR]. 
Section 9 investigates cycles in bases of lines of modular lattices $L$. It builds on point-splittings (Sec. 5) to retrieve relations among various invariants of $L$, e.g. a sharp lower bound for the number of join irreducibles of $L$. Given any BOL $\cal B$ of $L$ and any covering $(a,b)$ of $L$, it also pays to investigate the 'localization' $\B(a,b)$ of  $\cal B$ to $(a,b)$.

 Andr\'as Huhn [Hu]
 introduced the 2-distributive law, which is a wide-ranging generalization of (ordinary) 1-distributivity. It is obvious  that a modular and 2-distributive lattice doesn't contain (as interval the subspace lattice of) a projective plane, and Huhn managed to prove the converse. It turns out (Section 10) that the BOL's of 2-distributive modular lattices behave much better than in the general case. Namely, all BOLs $\B$ of $L$ are 'locally acyclic' (i.e. each localization $\B(a,b)$ of $\B$ is acyclic) iff $L$ is 2-disributive. 
 Here comes the result that kindled\footnote{Being a visiting Scholar at MIT (sponsored by the Swiss National Fonds)  Rota asked me to 'teach' him modular lattices on the blackboard in his office. Preferably in German, so that he could freshen up that too.} Rota's love of modular lattices: If $R$ is a suitably restricted ring (e.g. a skew-field) and $L$ is 2-distributive, then there is a lattice embedding of $L$ into the submodule lattice of a free $R$-module [HW]. As to other kinds of embeddings, the embeddability of arbitrary lattices into partition lattices is a popular theme in lattice theory soon after its inception. The author contributed [W3,W8] to this line of research by focusing on {\it cover-preserving (cp)} embeddings of {\it modular} lattices $L$ into partition lattices. It is easy to see that $L$ must be 2-distributive and without covering sublattices $M_4$, but this is not quite sufficient. Coming back to rings $R$,  if a modular lattice $L$ (of length $n$) is cp  embeddable into a partition lattice then, for each field $R=k$, $L$ is embeddable in $\La(k^n)$, the lattice of subspaces of $k^n$. Again, the converse fails.  As opposed to existence, Section 11 focuses on the {\it classification} of $k$-linear representations. Here we need to restrict locally acyclic (=2-distributive) lattices to certain globally acyclic lattices. Section 11 mainly  draws on [HW] and [W1].  The  axiomatization of bases of lines (Section 12), spearheaded in [HPR], simplifies a lot when attention is restricted (at least initially) to the 2-distributive case. Likely Section 12 will develop only with the help of co-authors.
  
Part D so far comprises two 'loose ends' but more may follow. Section 13 hopefully will present, the first time in English (!), deep results of Herrmann [H] and Wille [Wi1, Wi2]. In particular [Wi1] completely classifies those finite posets $P$ for which the freely generated modular lattice $FM(P)$ is finite. Section 14 tackles the 'practical side' of Wille's Theorem,
i.e. the compact (based on wildcards as in Sec. 8) representation of $FM(P)$. This is new material, even for the free {\it distributive} lattice $FD(P)$.

\bigskip
\centerline{\bf 2) Basic facts }
\bigskip
{\bf 2.1.} Some basic concepts (such as  dual lattice, direct products, congruences, etc) will not be defined in the sequel; their definitions can be found (most detailed)
 in [DP]; see also [G],[B],[D], which dig deeper.  For elements $x,y$ of a lattice $L$ we usually (not always) scorn the traditional notation
 $x\vee y$ and $x\wedge y$ for joins and meets and rather write 
  $x+y$ and $xy$
 respectively. In so doing, akin to  ordinary arithmetic, we can save brackets by adopting the convention that meets (=products) bind stronger than joins (=sums). Thus $x+yz$ means $x+(yz)$ and not $(x+y)z$.
 For $n\in \N$ we put $[n]:=\{1,2,...,n\}$.
 
A \underbar{chain} $C$ in a lattice is a totally ordered sublattice, i.e. for all $x,y\in C$ one either has $x\le y$ or $y\le x$. A lattice $L$ has \underbar{finite length (FL)} if
the cardinalities of all its chains  are bounded
by a universal constant. Notice that a FL lattice may well be infinite, such as $M_n$ defined below.

Recall that by convention (Sec. 1) all lattices are assumed to be of {\bf finite length (FL)}. We usually leave it to the reader to ponder in which scenarios this restriction could be weakened. It e.g. follows that for us each lattice $L$ has a {\it unit} (=largest element) which we denote by $1_L$ or simply $1$. Similarly $L$ has {\it zero} (=smallest element) which we write as $0_L$ or $0$.
Call $x\in L$  a \underbar{lower cover} of $y\in L$ (and $y$ an \underbar{upper cover} of $x$), written
$x\prec y$, if $x<y$ and there is no $z$ with $x<z<y$. For each $a<b$ there are maximal chains between $a$ and $b$, all of which necessarily of  type $a=x_0\prec x_1\prec\cdots\prec x_t=b$.
By definition the \underbar{length} of a $(t+1)$-element chain (maximal or not)  is $t$. The lattice $N_5$, the first in Figue 2.1, shows that not all maximal $a,b$-chains need to have the same length. In any  lattice $L$ we define $\delta(L)$ as the length of a {\it longest} maximal 0,1-chain. So $\delta(N_5)=3$.

 A {\it nonzero} element $p\in L$ is
\underbar{join-irreducible} if $x+y=p$ implies $x=p$ or $y=p$.
Each join irreducible
 $p$ has a unique lower cover $p_*\prec p$. Graphically speaking, $p$ is join irreducible iff in the diagram of $L$ only one line reaches $p$ from below. Denote by $J(L)$ the set
of  join irreducibles of $L$. Put

 $j(L):=|J(L)|\  and\  J(x):=\{p\in J(L)|\ p\le x\}$
 
  for all $x\in L$. Then $x=\sum J(x)$ for all $x\in L$ (see eg [DP,2.46]). Consequently for each $x<y$ the set $J(x,y):=J(y)\setminus J(x)$ is nonempty. According to Ex.2A it holds that

\begin{enumerate}
	\item [(1)] $\quad j(L)\ge \delta(L)$ in any finite length lattice $L$.
\end{enumerate}

The set $At(L)$ of all \underbar{atoms} $0\prec p$ is a subset of $J(L)$. A lattice is \underbar{atomistic} if $J(L)=At(L)$. For any $n\ge 3$ (including infinite cardinals) the lattice $M_n$ consists of 0,1, and $n$ atoms. The case $n=3$ yields the middle lattice in Figure 2.1. One defines 
\underbar{meet-irreducibles} and \underbar{co-atoms} in dual fashion. Actually, the \underbar{dual} $L^d$ of a lattice $L$ is obtained by turning the diagram of $L$ on its head, see [DP,1.19]. It then holds that the meet-irreducibles of $L$ match the join-irreducibles of $L^d$, and likewise for many other concepts.

An element $p\in L$ is
\underbar{join-prime} if $x+y\ge p$ implies $x\ge p$ or $y\ge p$, for all $x,y\in L$. It is immediate that join-prime $\Ra$ join-irreducible. Dually $m$ is \underbar{meet-prime} if $xy\le m$ implies $x\le m$ or $y\le m$, for all $x,y\in L$. It similarly follows that meet-prime $\Ra$ meet-irreducible.
Exercise 2B asks to determine the join-primes (and coupled meet-primes) in the three lattices of Figure 2.1.

\begin{center}
	\includegraphics[scale=0.9]{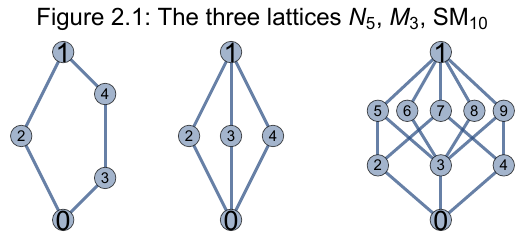} 
\end{center}

{\bf 2.2} Several types of functions between lattices $L,L'$ will be considered. A map $f:L\to L'$ is \underbar{monotone} if $a\le b\Ra f(a)\le f(b)$ for all $a,b\in L$. It is an \underbar{order 
	embedding} if $a\le b\Leftrightarrow f(a)\le f(b)$ for all $a,b\in L$. (Such $f$ is necessarily injective.)
For instance, let $L$ be the left lattice in Figure 2.2 and $L'$ be the powerset $\Pow([3])$ on the right. Defining $f$ by

$(2)\qquad 0\mapsto \es,\ a\mapsto \{1\},\ b\mapsto \{3\},\ c\mapsto \{1,2\},\ d\mapsto \{2,3\},\ e\mapsto \{1,2,3\}$

yields an order embedding. For instance $a\not\le d$ matches $\{1\}\not\s\{2,3\}$.
A \underbar{meet embedding} is an order embedding that additionally satisfies $f(ab)=f(a)f(b)$ for all $a,b\in L$.
A \underbar{join embedding} is an order embedding that additionally satisfies $f(a+b)=f(a)+f(b)$ for all $a,b\in L$. The order embedding $f$ defined in (2) is neither a meet embedding ($f(cd)=f(0)=\es\ne \{2\}=f(c)f(d)$), nor a join embedding ($f(a+b)=f(1)=\{1,2,3\}\ne \{1,3\}=f(a)+f(b)$).
However (Exercise 2C), each {\it bijective} order
embedding must be both a meet and a join embedding. 

\begin{center}
	\includegraphics[scale=0.93]{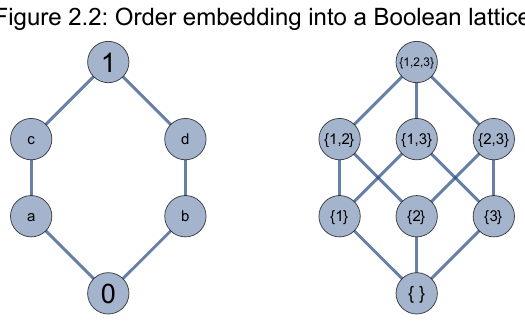} 
\end{center}

A \underbar{homomorphism} between lattices $L$ and $L'$ is any (e.g. constant) function $f:L\to L'$ satisfying $f(ab)=f(a)f(b)$
and $f(a+b)=f(a)+f(b)$ for all $a,b\in L$. Hence an injective homomorphism is the same thing\footnote{ {\it Non-injective} maps that preserve $\wedge$ (or $\vee$) are necessarily monotone but the extra property is seldom relevant.} as a simultaneous meet and join embedding. A surjective homomorphism is called an \underbar{epimorphism} and a bijective homomorphism is an \underbar{isomorphism}. Epimorphisms are linked to congruences, which will be the topic of Section 3.

{\bf 2.3} The \underbar{interval} determined
by elements $c<d$
of a lattice $L$ is the sublattice
\hfill\break $[c,d]:=\{x\in L:\ c\le x\le d\}$.    An interval
$[c,d]$ is \underbar{upper transposed} to the interval $[a,b]$,
written
$[a,b]\nearrow  [c,d]$, if
$a=bc$ and $d=b+c$. Dually $[c,d]$ is \underbar{lower transposed} to
$[a,b]$, written $[a,b]\searrow  [c,d]$, if $c=ad$ and $b=a+d$. For instance $[2,1]\searrow [0,4]$ in $N_5$,
and $[0,2]\nearrow [4,1]\searrow [0,3]$ in\footnote{Of course $[0,2]\nearrow [4,1]\searrow [0,3]$  abbreviates   $[0,2]\nearrow [4,1]$ {\it and}     $[4,1]\searrow [0,3]$. } $M_3$.
We may also use $\swarrow$ (equivalent to $\nearrow$) and $\nwarrow$ (equivalent to $\searrow$).
Being (upper or
lower) \underbar{transposed} is  denoted by $\sim$. The proof of the following is left as Exercise 2G.

{\bf Theorem 2.1: }{\it If  $L$ is a lattice and  $p\in L$ is join-prime, then there is a largest (not just maximal) element $m\in L$ with the property that $m\not\ge p$. 
	Hence $L=[p,1]\uplus [0,m]$.
	Moreover,  $m$ is meet-prime and $[p_*,p]\nearrow  [m,m^*]$.}

{\bf 2.4.} Before we come to particular types of lattices, a few words on three  tightly interwoven concepts are in order (more details can e.g. be found in [W7]):
 closure operators, closure systems, implicational bases. A \underbar{closure system} on a set $E$ is a family ${\cal F}$ of subsets stable under arbitrary intersections. Formally ${\cal F}\s\Pow(E)$  and $\bigcap{\cal G}\in {\cal F}$ for all ${\cal G}\s{\cal F}$. In particular $E=\bigcap\es\in {\cal F}$. 
 
  A \underbar{closure operator} is a function $cl:\Pow(E)\to \Pow(E)$ which is extensive, monotone, and idempotent. We  refer to the pair $(E,cl)$ as a \underbar{closure space}. The family $\La(E,cl)$ of all \underbar{closed} sets $X$ (i.e. $cl(X)=X$) is a closure system. Each closure system, ordered by set inclusion, is a lattice with meets and joins given by $X\wedge Y=X\cap Y$ and $X\vee Y=cl(X\cup Y)$ respectively.
  Conversely, each closure system ${\cal F}$ yields a closure operator by putting $cl(X):=\bigcap\{Y\in{\cal F}:\ X\s Y\}$. The two transitions are mutually inverse.
  
  The third concept amounts to pure Horn functions (a type of Boolean functions), but the following stripped-down formalism is more user-friendly. Call a pair of subsets $(A,B)\in {\cal P}(E)\times{\cal P}(E)$
   an \underbar{implication}. It is henceforth more handy to write $A\to B$ instead of $(A,B)$. We call $A$ the \underbar{premise} and $B$ the \underbar{conclusion} of the implication.
   Given a family $\Sigma=\{A_1\to B_1,\ldots, A_n\to B_n\}$  of implications on $E$,
  the set $X\s E$ is \underbar{$\Sigma$-closed} if for all $A_i\to B_i$ in $\Sigma$ it follows from $A_i\s X$ that $B_i\s X$. It is easy to see (Exercise 2H) that the family $\La(E,\Sigma)$ of all $\Sigma$-closed sets is a closure system. The coupled closure operator $cl_{\Sigma}(Y)$ is calculated by iteratively adding conclusions $B_i$ whenever the corresponding premise $A_i$ is contained in the superset of $Y$ obtained so far; see [W7] for concrete examples. Conversely, given any closure operator $cl:{\cal P}(E)\to {\cal P}(E)$, there always are families of implications $\Sigma$ such that $cl=cl_{\Sigma}$. We then call $\Sigma$ an \underbar{implicational base} of $cl$, or of the coupled closure system $\La(E,cl)$.

\medskip

{\bf 2.5.}
A lattice $D$ is \underbar{distributive} if $x(y+z)=xy+xz$ for all
 $x,y,z\in D$. Then automatically the dual identity holds: $x+yz=(x+y)(x+z)$ for all $x,y,z\in L$. For instance (Ex. 2I), each chain  is a distributive lattice. The 2-element chain $D_2$ will come up frequently. If $(J,\le)$ is a poset, then the collection $\La(J,\le)$ of all its order ideals is easily seen to be closed under intersections and unions. Hence  $\La(J,\le)$ 
 is a sublattice of the distributive powerset lattice $\Pow (J)$, and so itself is distributive.
 
  The converse also holds (Birkhoff's Theorem, [DP,p.116]): Each (FL) distributive  $L$ is actually finite and the map $x\mapsto J(x)$ is an isomorphism from $L$ to  $\La(J(L),\le)$.
  As a consequence, each element $x$ of a distributive lattice has a {\it unique} irredundant representation by join-irreducibles (i.e. by the set of maximal elements of $J(x)$). By duality also the irredundant representations by meet-irreducibles are unique.
   Replacing each set $J(x)\s J(L)$ by its characteristic 0,1-vector, and using that all join irreducibles are join-prime, one can view (Ex.2O) the map $x\mapsto J(x)$ as a subdirect product representation $L\to (D_2)^n$. 
 
  A lattice $L$  is \underbar{complemented} if each $x\in L$ has a \underbar{complement} $c$, i.e. $x+c=1$ and $xc=0$. If  $L$ is  \underbar{Boolean}, i.e. distributive and  complemented, then $(J(L),\le)$ is unordered (i.e. an \underbar{antichain}), and so $\La(J(L),\le)$ becomes $\Pow(J(L))$.
 
\bigskip

{\bf 2.6} If an interval $[c,d]$ has cardinality two, i.e. $c\prec d$, then we call\footnote{Other authors speak of  prime intervals or of prime quotients.}
it a \underbar{covering}. A  lattice $L$ is \underbar{upper semimodular} if upper transposes of coverings  are again coverings . For instance, the rightmost lattice in Figure 2.1, call it $SM_{10}$, is upper semimodular, e.g. the covering  $[0,2]$ transposes up to the covering  $[8,1]$. Notice that coverings  can have {\it lower} transposes which are not coverings, e.g. $[7,1]\searrow  [0,6]$ in $SM_{10}$. An upper semimodular lattice in particular has the property that from $a\prec b,\ a\prec c$ it follows that $b\prec b+c,\ c\prec b+c$. In fact this property is equivalent to our definition.
Evidently each {\it interval} sublattice of an upper semimodular lattice is again upper semimodular, but inheritance of upper semimodularity is not\footnote{In fact, a famous Theorem of Dilworth states that {\it each} finite lattice occurs as sublattice in a suitable semimodular lattice. A short proof can be found in [W3].} guaranted for arbitrary sublattices.

Each upper semimodular lattice has a \underbar{rank function}
$\delta:L\to \N$, i.e. satisfying $\delta(0)=0$ and $\ x\prec y\Rightarrow\delta(y)=
\delta(x)+1$. The existence of $\delta$ is due to the fact that all maximal chains $C$ between elements $a<b$ 
have the same length (which e.g. fails in $N_5$). In particular, the \underbar{length of $L$}, defined as $\delta(L):=\delta(1)$, is the length of any maximal 0-1-chain. Furthermore, it holds that

 $(3)\qquad \delta(x)-\delta(xy)\ge\delta(x+y)-\delta(y)$,

i.e. switching from an interval to an upper transposed interval cannot increase the length\footnote{Many textbooks state this inequality in equivalent (yet pedagogically inferior) fashion as $ \delta(xy)+\delta(x+y)\le \delta(x)+\delta(y)$.}. 

{\bf 2.6.1} A lattice $L$ is \underbar{lower semimodular} if lower transposes of coverings are coverings. Not surprisingly, in this case $\ge$ in (3) becomes $\le$.
Put another way, $L$ is lower semimodular iff $L^d$ is upper semimodular. Henceforth upper semimodular will be more important than lower semimodular, and so we convene that \underbar{semimodular} always abbreviates 'upper semimodular'.
The  lattice  $SM_{10}$ in Figure 2.1 is semimodular but  not lower semimodular. The lattice $N_5$ is neither semimodular nor lower semimodular. We leave it as an exercise to prove the following fact [Lemma 1 in W3].

\begin{enumerate}
	\item[(4)] Let $f:L\to L'$ be a cover preserving order between lattices. If $L$ semimodular, then $f$ preserves joins. Dually, if $L$ lower semimodular, then $f$ preserves meets. 
\end{enumerate}

The example in Section 2.2 shows that the (upper resp. lower) semimodularity of $L$ in (4) is essential.

{\bf 2.7} By definition a  lattice $L$ is \underbar{modular} if it is both upper and lower semimodular. Consequently (3) improves to the \underbar{rank formula} 

\begin{enumerate}
	\item[(5)] $\delta(x)-\delta(xy)=\delta(x+y)-\delta(y).$ 
\end{enumerate} 

A handy consequence (Exercise 2Q) is the following:

\begin{enumerate}
	 \item[(6)] Suppose in a modular lattice $b_i\le a_i$ and $d_i:=\delta(a_i)-\delta(b_i)\ (i=1,2)$. \\Then $\delta(a_1+a_2)-\delta(b_1+b_2)\le d_1+d_2$.
\end{enumerate}

Transposed intervals not only have the same length, they are in fact isomorphic:

\begin{enumerate}
	\item[(7)] {\it Dedekind's transposition principle:  Whenever  $[a,b]\nearrow  [c,d]$ takes place in a modular lattice, the
		map $x\mapsto x+c$ is a
		lattice isomorphism between $[a,b]$ and $[c,d]$ with inverse $y\mapsto
		yb$. } 	
\end{enumerate}

Sometimes useful is this consequence [FRS, p.361] of (7) which generalizes the well-known (and easy) case that $p,q\in At(L)$:

\begin{enumerate}
	\item[(8)] {\it If $L$ is modular and $x\in L,\ p,q\in J(L)$ are such that $p\le x+q$ but $p\not\le x+q_*$,\\ then $x+p=x+q$. } 	
\end{enumerate}

To prove it,  $y\mapsto x+y$ is an isomorphism $[xq,q]\to [x,x+q]$ between {\it non-singleton} intervals (why?). Because therefore $q_*$ is the unique lower cover of $q$ in $[xq,q]$, we have $x+q_*$ as the unique lower cover of $x+q$ in
$[x,x+q]$. By assumption $x+p\in [x,x+q]$, but $x+p\not\in [x,x+q_*]$. This forces $x+p=x+q$ and proves (8).

The proof of the  Transposition Principle (7) hinges on the fact that 

\begin{enumerate}
	\item[(9)] $(x+y)z=x+yz$ for all $x,y,z\in L$ 
with $x\le z$.
\end{enumerate}

Condition (9) is actually the usual definition of {\it modularity}, but we hasten to add that in our FL framework (9) is equivalent to 'semimodular plus lower semimodular'. It is easy to see that the conditioned (by $x\le z$) identity appearing in (9) is equivalent to the unconditioned identity $(xz+y)z=xz+xy$. This observation e.g. implies  that the family of all modular lattices is a \underbar{variety}, i.e. closed under taking sublattices, direct products, and epimorphic images. Yet varieties will concern us only in Part D of  this text.
As shown in Exercise 2M:

\begin{enumerate}
	\item[(10)] If $(a,b)$ is a covering in a modular lattice then $J(a,b)$ is an antichain.
\end{enumerate}

(This property e.g. fails in $N_5$.) Subspace lattices of vector spaces, and their sublattices, are among the most important examples of modular lattices.
Also each distributive lattice is modular because $(xz+y)z=xzz+yz=xz+xy$. As is well known [DP,4.10], a modular lattice  is
distributive  iff it does not contain a sublattice $M_3$.
We will see soon that due to our FL convention  $M_3$ can be found as a {\it covering} (= length two) sublattice.
 
A lattice is \underbar{sectionally complemented} if all its intervals are complemented.

{\bf Theorem 2.2: }{\it For each modular lattice $L$ these properties are equivalent:
\begin{enumerate}
	\item[(a)] It is complemented.
	\item[(b)] It is sectionally complemented.
	\item[(c)] It is atomistic.
	\item[(d)] Its unit is a sum of atoms.
	\item[(e)] It is co-atomistic.
	\item[(f)] Its zero is a meet of co-atoms.
	
\end{enumerate}
}

{\it Proof.} As to $(a)\Ra (b)$, in order to find a complement of $x\in [a,b]$ within this interval, let $x'$ be a complement within $L$ and put $y:=a+x'b\in [a,b]$. Using (8) we conclude that $x+y=x+x'b=(x+x')b=1b=b$ and $yx=(a+x'b)x=a+x'bx=a+0=a$.
As to $(b)\Ra (c)$, if any lattice has a join irreducible $p$ which is not an atom, then $p_*\ne 0$ has no complement in the interval $[0,p]$. Implication  $(c)\Ra (d)$ is trivial,
and $(d)\Ra (a)$ is covered by Exercise 2L. 
 Since the dual of a modular lattice is modular, we have $(c)\LRa (e)$ and $(d)\LRa (f)$. $\square$

 For instance, since the modular lattice in Figure 4.1 is such that 1 is the join of atoms, it must be complemented.
 Furthermore, we will e.g. exploit the following: In each modular lattice $L$ (complemented or not), when $x_0$ is the meet of all lower covers of $x\in L$, then $[x_0,x]$ is complemented (e.g. Boolean when $L$ was distributive). The fine structure of complemented modular lattices will occupy us in Section 4.  It is interesting to pinpoint (Ex.2L) what survives from Theorem 2.2 when $L$ is merely semimodular.

\bigskip
{\bf 2.8 Exercises.}
\bigskip

{\begin{enumerate}
	\item[\bf Ex. 2A.] Show that always $j(L)\ge \delta(L)$.
	
	{\it Solution of Exercise 2A.} Let $0\prec x_1\prec x_2\prec\cdots\prec x_n=1$ be a longest chain in $L$ (so $n=\delta(L)$). For all $1\le i\le n$ pick any $p_i$ in $J(x_{i-1},x_i)$ (which is nonempty by 2.1). The $p_i$'s are distinct because for all $i<j$ it holds that $x_{j-1}+p_j=x_j$ but $x_{j-1}+p_i=x_{j-1}$. Hence $j(L)\ge n$.
		
		\item[\bf Ex. 2B.] Determine all join-prime elements $p$ and their coupled meet-primes $m$ for the lattices in Figure 2.1.
		
		{\it Solution of Exercise 2B.} In $N_5$ the asked pairs $(p,m)$ are (2,4) and (3,2), in $SM_{10}$ it is (3,7), and $M_3$ has no join-prime elements.

		\item[\bf Ex. 2C. ]Show that every bijective\footnote{Strictly speaking, talking of an 'embedding' is inappropriate here since $f(L)$ is no smaller than $L'$. More precise is it to talk of a bijective map $f$ that satisfies  $a\le b\LRa f(a)\le f(b)$ for all $a,b\in L$. But then again, adopting this pedantry (or pedanticism, pedantism?!) would also force us to trim the definition of 'order embedding' in 2.2 by adding 'non-bijective'. } order embedding $f:L\to L'$ is a lattice isomorphism.
		
		{\it Solution to Exercise 2C.} Let $f:L\to L'$ be such that $a\le b\LRa f(a)\le f(b)$ for all $a,b\in L$. Since $f$ is bijective, it has an inverse $f^{-1}$ which clearly satisfies $\alpha\le \beta\LRa f^{-1}(\alpha)\le f^{-1}(\beta)$ for all $\alpha,\beta\in L'$. From  $x,y\le x+y$ follows $f(x),f(y)\le f(x+y)$, hence $f(x)+f(y)\le f(x+y)$. We now show that
		
		$$x+y\quad {{(\&)}\atop{\le}}\quad f^{-1}(f(x)+f(y))\quad {{(\#)}\atop{\le}}\quad x+y.$$
		
		This immediately gives $f(x+y)=f(x)+f(y)$ (and similarly $f(xy)=f(x)f(y)$). As to $(\&)$, 
		from $f(x),f(y)\le f(x)+f(y)$ follows $x,y\le f^{-1}(f(x)+f(y))$, hence $x+y\le f^{-1}(f(x)+f(y))$. As to $(\#)$, this follows from $f(x)+f(y)\le f(x+y)$.
		
		\item[\bf Ex. 2D.] Show that $f(x):=J(x)$ is a meet embedding $L\to\Pow(J(L))$ for each lattice $L$. 
		
		{\it Solution to Exercise 2D.} If $J(x)\not\s J(y)$, there is $p\in J(x)\setminus J(y)$, and so $x\not\le y$. Hence $f$ is an order embedding. As to meet-preservation, clearly $J(x\wedge y)\s J(x)\cap J(y)$. Conversely, pick $p\in J(x)\cap J(y)$. Then $p$ is a common lower bound of $x$ and $y$. But  the largest common lower bound is $x\wedge y$, and so $p\le x\wedge y$.
		
		\item[\bf Ex. 2E] Let $\varphi:L\to L'$ be a homomorphism. Show that  $\varphi(x)<\varphi(y)$ does not imply $x<y$. (Hint: Take $L:=N_5$). However, show that $\varphi(x)<\varphi(y)$ {\it does} imply
		$\varphi(x)=\varphi(z)$ for some $z<x$.
		
		{\it Solution of Exercise 2E.} Define $\varphi:N_5\to D_2$ by $\varphi(0)=\varphi(2)=0$ and $\varphi(a)=1$ otherwise. Then $\varphi(2)<\varphi(3)$ but $2\not<3$. Consider now any homomorphism $\varphi:L\to L'$ with $\varphi(x)<\varphi(y)$. Then $\varphi(xy)=\varphi(x)\varphi(y)=\varphi(x)\neq\varphi(y)$, hence $xy\neq y$, hence $z:=xy<y$ satisfies $\varphi(z)=\varphi(x)$.
		
		\item[\bf Ex. 2F.] Show that from $[a,b]\nearrow [c,d]\nearrow  [e,f]$ follows $[a,b]\nearrow  [e,f]$.
		
		{\it Solution of Exercise 2F.} We need to show $be= a$ and $b+e= f$. As to the first identity, $be\ge bc=a$, and
		$be\le de=c$ implies $be\le bc=a$. Hence $be=a$. Similarly one shows $b+e= f$.
		
		\item[\bf Ex. 2G.] Prove Theorem 2.1.

		{\it Solution to Exercise 2G.}  Let $p$ be join-prime. Using FL there are elements $m,m',\ldots$ maximal with the property that $m\not\ge p,\ m'\not\ge p$, and so forth. Hence $m+m'\not\ge p$ by definition of join-prime. By maximality 
		$m=m'$ is {\it the largest} element $\not\ge p$ in $L$. In other words,
		$L=[p,1]\uplus [0,m]$. In order to show that $m$ is meet-prime (whence a fortiori meet irreducible), take any $x,y\not\le m$. Then $x,y\in [p,1]$, hence $xy\in [p,1]$, hence $xy\not\le m$.
		
		From $p\le m^*,\ p\not\le m,\ m\prec m^*$ follows $p+m=m^*$. Similarly $pm=p_*$, and so $[p_*,p]\nearrow  [m,m^*]$.

		\item[\bf Ex. 2H.] Given a family $\Sigma$ of implications on $E$, show that the system $\La(E,\Sigma)$ of all $\Sigma$-closed sets is closed under intersections. If additionally all implications in $\Sigma$ have singleton premises, show that $\La(E,\Sigma)$ is closed under unions as well.
		
		{\it Solution to Exercise 2H.} Let $X,Y\s E$ be $\Sigma$-closed and take any $A\to B$ in $\Sigma$.
		If $A\not\s X\cap Y$, there is nothing to show. If $A\s X\cap Y$, then $A\s X$ and $A\s Y$, hence $B\s X$ and $B\s Y$, hence $B\s X\cap Y$. Additionally we now assume that all implications in $\Sigma$ have singleton premises. Again, let $X,Y\s E$ be $\Sigma$-closed and take any $A\to B$ in $\Sigma$
		with $A\s X\cup Y$. Generally this neither implies $A\s X$ nor $A\s Y$, but here it does since $A=\{a\}$. If say $A\s X$, then $B\s X\s X\cup Y$.

		\item[\bf Ex. 2I] Show that each chain is a distributive lattice.

		{\it Solution to Exercise 2I.} For all $x,y$ in a chain one has $x\vee y=max(x,y),\ x\wedge y=min(x,y)$.
		If say $z<x<y$, then $x\wedge (y\vee z)=x\wedge y=x=x\vee z=(x\wedge y)\vee (x\wedge z)$. For the other five orderings of $x,y,z$ one proceeds likewise.

		\item[\bf Ex. 2J] (From [DP,p.125,Ex.5.5])  Using Ex. 2I, and the fact that distributivity carries over to direct products,  show that the poset of all divisors of a natural number $n$ (partially ordered by $x\le y:\LRa y=xz$) is a distributive lattice.
		
		\item[\bf Ex. 2K] Show that in a finite distributive lattice all join irreducibles are join-prime.
		
		{\it Solution to Exercise 2K.} Let $p\in J(L)$ and let $x,y\in L$ be such that $p\le x+y$. Then, using distributivity, $p=p(x+y)=px+py$. Since $p$ is join irreducible, $px=p$ or $py=p$, hence $x\ge p$ or $y\ge p$.

		\item[\bf Ex. 2L] This exercise appeals to Theorem 2.2. Let $L$ be a semimodular lattice whose unit $1$ is a sum of atoms. Show that:
		\begin{enumerate}
			\item[(i)] $L$ needs not be  atomistic.	(Hint: Consider $SM_{10}$ in Figure 2.1.)
			\item[(ii)] $L$ is co-atomistic.
			\item[(iii)] $L$ is complemented.	
		\end{enumerate}
		
		{\it Solution to Exercise 2L.} 
		
		As to (ii), we need to derive a contradiction from assuming there is a meet irreducible $m$ with $m^*\ne 1$. Since $1$ is a sum of atoms, there is an atom $p$ with $p\not\le m^*$.
		A fortiori $p\not\le m$, and so semimodularity implies $m\prec m+p$. But then $m+p\ne m^*$ is a second upper cover of $m$, contradicting meet irreducibility.
		
		As to (iii), to fix ideas suppose $d(L)=8$ and $d(x)=5$. We calculate a complement $x'$ of $x$ as follows. Since $1$ is a sum of atoms, there is $0\prec p_1$ with $p_1\not\le x$, and so $d(x+p_1)=6$ by semimodularity.
		Similarly there are atoms $p_2\not\le x+p_1$ and $p_3\not\le x+p_1+p_2$. It follows that $d(x+p_1+p_2+p_3)=8$.
		Putting $x':=p_1+p_2+p_3$ we have $d(x')\le 3$ and $x+x'=1$, as well as $d(xx')\le d(x)+d(x')-d(x+x')\le 5+3-8=0$. Hence $xx'=0$ (and $d(x')=3$).
		
	\item[\bf Ex.2M.] Show that in a  modular lattice $L$ each join irreducible $p$ is 'weakly' join-prime in the sense  that from $p\le p_*+x$ follows $x\ge p$. Furthermore proof claim (10).
	
	{\it Solution to Exercise 2M.} Suppose that $x+p_* \ge p$. Then $p=(x+p_*)p=p_*+xp$. From $p\in J(L)$ follows $xp=p$, i.e. $x\ge p$. As to (10), suppose there were $p,q\in J(a,b)$ with $p<q$. Then $p\le q_*$ implies $q\le b=a+p\le a+q_*$, contradicting the weak primality of $q$.

	\item[\bf Ex.2N.] For each vector space $V$ (over any field and possibly of infinite dimension) let $\La(V)$ be the closure system of all  subspaces. Convince yourself that $X\vee Y$ is the usual sum $X+Y$ of subspaces. Show that $\La(V)$ is modular by applying (9) rather than upper plus lower semimodularity.
	
	{\it Solution to Exercise 2N.} Each subspace containing $X$ and $Y$ must contain $X+Y$. Since $X+Y$ is itself a subspace, it coincides with the join (:=least common upper bound) $X\vee Y$.
	
	In {\it any} lattice it follows from $X\le Z$ that $(X\vee Y)\wedge Z\ge X\vee (Y\wedge Z)$ (check). It hence suffices to verify that in the situation of vector subspaces the converse inequality $(X+Y)\cap Z\s X+(Y\cap Z)$ holds. Indeed, every element on the left is of the form $x+y$, where $x\in X,\ y\in Y,\ x+y\in Z$. From $x\in X\s Z$ follows $y=(x+y)-x\in Z$, and so $x+y\in X+(Y\cap Z)$.
	
	\item[\bf Ex.2O.] Let $L$ be a lattice all of whose join-irreducibles $p_1,\ldots , p_s$ are join-prime. For each $x\in L$ define $f(x)$ as the bitstring of length $s$ whose $i$-th component is 1 iff $p_i\le x$. Show that $f$ provides a subdirect representation $L\to (D_2)^s$. Conclude that $L$ must be distributive. Hint: Use exercises 2D.
	
	\item[\bf Ex.2P.]  Show that distributivity implies $(a+b)(a+c)(b+c)=ab+ac+bc$, and modularity implies
	$(a+b)(a+c)(b+c)=b(a+c)+a(b+c)$.
	
	\item[\bf Ex. 2Q] Prove claim (6).
	
		{\it Solution to Exercise 2Q.} We have
		
		 $\delta(a_1+a_2)-\delta(b_1+b_2)=\delta(a_1)+\delta(a_2)-\delta(a_1a_2)-( \delta(b_1)+\delta(b_2)-\delta(b_1b_2))=d_1+d_2-\delta(a_1a_2)+\delta(b_1b_2)\le d_1+d_2$, viewing that $b_1b_2\le a_1a_2$ implies $\delta(b_1b_2)-\delta(a_1a_2)\le 0$.

\end{enumerate}
\bigskip

\centerline{\bf 3.  Lattice congruences via coverings }
\medskip
 Recall from  universal algebra [KNT] that a \underbar{congruence (relation)}
$\theta\s A\times A$ of a (universal)
algebra $A$ yields a factor algebra $A/\theta$ and that the set of all congruences
 forms a lattice $Con(A)$. Here $\theta\wedge \theta'=\theta\cap
\theta'$ and $\theta\vee\theta'$ is the smallest congruence containing $\theta\cup
\theta'$. The $0$ of $Con(A)$ is the identity relation $\{(x,x):\ x\in A\}$ and $1\in Con(A)$ is the universal relation $A\times A$.
If $Con(A)=\{0,1\}$, then $A$ is  \underbar{congruence-simple} (in brief: \underbar{c-simple}). 
By the second isomorphism theorem, $A/\theta$ is c-simple iff $\theta$ is a co-atom in $Con(A)$.
We will be mostly concerned with lattices $A=L$.

\medskip

{\bf 3.1.}  Here comes an important example. Suppose that $p\in L$ is join-prime with coupled meet-prime $m$. Recall from Theorem 2.1 that $L=[p,1]\uplus [0,m]$. Define the map $f:L\to \{0,1\}$ by $f(x):=1$ for $x\in [p,1]$ and $f(x):=0$ for $x\in [0,m]$. Exercise 3B shows that $f$ is an epimorphism whose kernel $\theta\in Con(L)$ has the congruence classes $[p,1]$ and $[0,m]$.
 Conversely, let $f:L\to D_2$ be any epimorphism with kernel $\theta$. Then the smallest element of $f^{-1}(1)$ is a join-prime element of $L$.

{\bf 3.2} 
Let $\theta\in Con(L)$. Then the \underbar{blocks} (=equivalence classes) $S$ of $\theta$ are easily seen to have these properties:
\begin{enumerate}
	\item[(1)] Each block is an interval\footnote{If $L$ is not assumed to be FL (as it isn't in [DP]) then 'interval' needs to be replaced by 'convex sublattice'.} of $L$.
	
	\item[(2)] The blocks $S$ are closed under transposition in this sense. If $[a,b]\s S$ and $c,d\in L$ are such that $[a,b]\sim [c,d]$, then $[c,d]$ is also contained in a block.
\end{enumerate}

Conversely [DP p.136], if  an equivalence relation $\theta$ on $L$  satisfies (1) and (2), then $\theta$ is a congruence. Here comes another characterization of congruences. Let $\theta$ be a reflexive binary relation on a lattice $L$. According to [KNT p.90] $\theta$ is a congruence on $L$ iff for all $a,b,c\in L$:
\begin{enumerate}
	\item[(3)] $a\theta b$ iff $(ab)\theta(a+b)$.
	\item[(4)] If $a\le b\le c$ and $a\theta b$ and $b\theta c$, then $a\theta c$.
	\item[(5)] If $a\le b$ and $a\theta b$, then $(ac)\theta(bc)$ and $(a+c)\theta(b+c)$ .	
\end{enumerate}

{\bf 3.3} Recall from Section 2.6 that intervals of cardinality two are called coverings. From now on we denote intervals $[a,b]$ that are coverings by $(a,b)$; then $(a,b)$ can conveniently be viewed as member\footnote{For arbitrary ordered pairs we use other letters, say (x,y).} of $L\times L$.
 Let $Cov(L)$  be the set of all coverings  of $L$. Furthermore, for $\theta\in Con(L)$ let $Cov(\theta)$ be the set of of all $(a,b)\in Cov(L)$ with $(a,b)\in\theta$.
 In the sequel we strive to reduce congruences of $L$ to the behaviour of $Cov(L)$. For starters, if
  $a<b$, there always are \underbar{chains of coverings } from $a$ to $b$ in the sense that $a=x_0\prec x_1\prec\cdots\prec x_n=b$. 
 Further it is clear from $(3),(4),(5)$ that for all $a,b\in L$ it holds that $a\theta b$ iff there is a chain of coverings, all of which belonging to $Cov(\theta)$, from $ab$ to $a+b$.
It follows that each $\theta\in Con(L)$ is uniquely determined by $Cov(\theta)$. In particular, $\theta\neq 0\Ra Cov(\theta)\neq\es$. Here comes an inherent characterizations of the sets of type $Cov(\theta)$.

{\bf Lemma 3.1: }{\it For each  lattice $L$ and each (possibly empty) set ${\cal Q}\s Cov(L)$ the following are equivalent:
\begin{enumerate}
	\item[(a)]  ${\cal Q}$ is \underbar{closed under subtransposition}, i.e. for all
	 $(a,b)\in {\cal Q}$ and all intervals $[c,d]$ with $(a,b)\sim [c,d]$ it follows that
	$(c',d')\in {\cal Q}$ for all $c\le c'\prec d'\le d$.
	\item[(b)] There is $\theta\in Con(L)$ with ${\cal Q}=Cov(\theta)$.
\end{enumerate} } 

{\it Proof.} The implication $(b)\Ra (a)$ holds because of (1) and (2). Assume now that ${\cal Q}\s Cov(L)$ satisfies (a). Our candidate $\theta\in Con(L)$ with ${\cal Q}= Cov(\theta)$ is the binary relation $\theta$ defined by $(a,b)\in\theta$ iff  there is a chain of coverings  from ${\cal Q}$ stretching from $ab$ to $a+b$ (for $a=b$ that chain is empty). It suffices to show that $\theta$ satisfies properties (3),(4),(5). Properties (3),(4) are immediate and (5) is left as exercise 3E. $\square$

{\bf Theorem 3.2: }{\it The congruence lattice $Con(L)$ is distributive for each lattice $L$.}

{\it Proof.} This is a well-known result of Funayama and Nakayama [DP,6.19]. The proof can be trimmed for our FL  lattices $L$. Namely, as previously noticed, then $\theta\in Con(L)$ is uniquely determined by $Cov(\theta)$, and so
  for all $\theta\in Con(L)$ it holds that
 $\theta\le\theta'$ iff $Cov(\theta)\s Cov(\theta')$. Hence the distributivity of $Con(L)$ will ensue if we  show that the family (=set system) $\{Cov(\theta):\ \theta\in Con(L)\}$ is stable under taking unions and intersections and whence, akin to $\La(J,\le)$ in 2.5, is a necessarily distributive sublattice of the powerset lattice of $Cov(L)$. By Lemma 3.1 it suffices to show that the family of all sets $\cal Q$ closed under subtransposition is stable under unions and intersections. Since the property that ${\cal Q}\s Cov(L)$ is closed w.r.to subtransposition readily translates to being $\Sigma$-closed w.r.to a certain family $\Sigma$ of implications with singleton premises, the claim follows from Exercise 2H.$\square$.
 
 {\bf 3.3.1.} Recall\footnote{For readers not familiar with the upcoming concepts we recommend [KNT].}  that a universal algebra $B$ is \underbar{subdirectly irreducible} if $Con(B)$ has a unique atom. Consequently $B=A/\theta$ is subdirectly irreducible iff $\theta\in Con(A)$ is meet-irreducible. A universal algebras $A$ may have several subdirect representations into subdirectly irreducible factors. Let now $A=L$ be a lattice. Since in  distributive lattices $D$ (such as $D=Con(L)$) the zero has a unique irredundant representation as meet of meet-irreducibles, $L$ has a unique representation $L\simeq L'\s L/\theta_1\times\cdots\times L/\theta_t$  with subdirectly irreducible factors $L/\theta_i$.
 What if $Con(L)$ is not just distributive but Boolean? Then $\theta_1,...,\theta_t$ above must be the co-atoms of $Con(L)$, and all $L/\theta_i$ must even be c-simple. This foreshadows Theorem 3.3.
 
 {\bf 3.4.} Let $\approx$ be the transitive closure of $\sim$ (transposition). Thus $\approx$ is an equivalence relation on the set $Int(L)$ of all nontrivial intervals of a lattice. Intervals $[a,b]\approx [c,d]$ are called \underbar{projective}. Unfortunately $[a,b]$ and $[c,d]$ need not be isomorphic, not even of the same length. In particular, if $(a,b)$ and $(c,d)$ are coverings such that $(a,b)\approx (c,d)$, then the {\it intermediate} intervals in a \underbar{witnessing sequence} of transpositions need not be coverings. This leads to the next definition.
 
  We   call $(a,b),(c,d)\in Cov(L)$  {\it homogeneous-projective} if all intermediate intervals in some\footnote{It is allowed that other witnessing sequences have intermediate intervals which are not coverings, see Exercise 3F} witnessing sequence $(a,b)\sim [u,v]\sim\cdots \sim (c,d)$ are  also coverings.
 By definition the\\ \underbar{hom-proj classes (of coverings)}  are the connected components of the graph with vertex set $Cov(L)$, where vertices are adjacent iff they are transposed.  For each  lattice we put
 \begin{enumerate}
 	\item[(6)] $\quad s(L)\ :=\ $ {\it the number of hom-proj classes constituting $Cov(L)$.}
 \end{enumerate}

Recall that $j(L):=|J(L)|$ (which may be infinite, say for $L=M_{\aleph_3}$). It is easy to see (Ex.3A) that

\begin{enumerate}
	\item[(7)] {\it $\quad j(L)\ge s(L)$ in each finite length lattice $L$.} 
\end{enumerate}
 
 For instance, $N_5$ has the hom-proj  classes of coverings
 $\{[4,1],[0,2]\},\ \{[2,1],[0,3]\},\ \{[3,4]\},$
  and so $j(N_5)=s(N_5)=3$. 
   We empasize that a hom-proj  class of coverings need not be closed  under sub-transposition (thus $[3,4]\not\in \{[4,1],[0,2]\}$ in $N_5$).
 
  However, let ${\cal Q}\s Cov(L)$ be closed under subtransposition.
  Then it follows from $(a,b),(c,d)\in Cov(L)$ and $(a,b)\in{\cal Q}$ and $(a,b)\sim (c,d)$ that $(c,d)\in {\cal Q}$. Therefore $\cal Q$ is a certain disjoint union of hom-proj classes. Hence the number of subtransposition-closed subsets of $Cov(L)$ is $\le 2^{s(L)}$. In view of Lemma 3.1 we conclude:
  
  \begin{enumerate}
  	\item[(8)] {\it For any lattice $L$ it holds that $|Con(L)|\le 2^{s(L)}$.} 
  \end{enumerate}

 {\bf Theorem 3.3: }{\it If $L$ is a modular lattice of finite length, then $Con(L)$ is a Boolean lattice
 	of length $s(L)$. Furthermore, if the lattices $L/\theta_i$ are the subdirectly irreducible factors of $L$ 
 	(i.e. the $\theta_i$'s are the co-atoms of $Con(L)$), then $\delta(L)=\sum_{i=1}^{s(L)}\delta(L/\theta_i)$. Finally, for each projectivity class ${\cal P}$ of coverings, each maximal 0-1-chain has the {\rm same} number of coverings  belonging to ${\cal P}$. }
 
 {\it Proof.} In a modular lattice coverings can only transpose to coverings. In other words, "closed under subtransposition" amounts to "closed under transposition". Because therefore each hom-proj class\footnote{One can succinctly put it this way. In a modular lattice the hom-proj classes of coverings are exactly the projectivity classes of intervals containing at least one covering.} is subtransposition-closed, {\it every} union of hom-proj classes is subtransposition-closed as well.
  It follows  that $Con(L)$ is a Boolean lattice of cardinality $2^{s(L)}$. 
 (As to the relevance of Boolean congruence lattices, see 3.3.1.)
 
As to the other claims, fix any maximal 0-1-chain $C$ in $L$. Let $d_i$ be
the number of coverings  in $C$ that belong to the $i$th projectivity class ${\cal P}_i$ of coverings . Trivially $\sum_{i=1}^s d_i=\delta(L)$. Let ${\cal Q}$ be the union of all projectivity classes of coverings  {\it distinct} from a fixed  ${\cal P}_i$. Then ${\cal Q}=Cov(\theta_i)$
for a unique co-atom $\theta_i=\theta({\cal P}_i)$ of $Con(L)$ which satisfies $\delta(L/\theta_i)=d_i$. (Thus $\theta_i$ {\it separates} exactly the coverings in ${\cal P}_i$.)
Therefore $d_i$ is  independent from the particular choice of $C$. $\square$

{\bf 3.4.2} Akin to $\delta(L)$ in Theorem 3.3, also $j(L)$  can be reduced to the subdirectly irreducible factors. Here come the details. Trivially each covering of type $(p_*,p)$ (so $p\in J(L)$) belongs to some projectivity class of coverings. Slightly more subtle (Ex.3A), each projectivity class of coverings $(c,d)$ contains at least one covering of type $(p_*,p)$. 

 By abuse of language we say that join irreducibles $p,q$ are \underbar{projective} if $[p_*,p]\approx [q_*,q]$. And we say that $p$ \underbar{belongs} to a projectivity class ${\cal P}$ of coverings  if $[p_*,p]\in {\cal P}$. To summarize:

\begin{enumerate}
	\item[(9)] {\it  Let $L$ be   modular and let ${\cal P}_1,\ldots,{\cal P}_s$ be its $s=s(L)$ many projectivity classes of coverings. Each set $J({\cal P}_i)$ of join irreducibles belonging to ${\cal P}_i$ is nonempty and\\ $J(L)=J({\cal P}_1)\uplus\cdots\uplus J({\cal P}_s)$. }
\end{enumerate}

Here comes the connection to the c-simple factors of $L$.

\begin{enumerate}
	\item [(10)] {\it  Let $\theta_i\in Con(L)$ be the co-atom that separates exactly the coverings in ${\cal P}_i$, and let $\varphi_i: L\to L/\theta_i$ be the canonical epimorphism. Then $\varphi_i$ maps  $J({\cal P}_i)$ bijectively upon $J(L/\theta_i)$, which implies  $ j(L)=\sum_{i=1}^{s(L)} j(L/\theta_i)$.   }
\end{enumerate}

{\it Proof of (10).} Introducing new notation it holds by the above that\\
 $J({\cal P}_i)=\{p\in L:\ \varphi_i(p_*)<\varphi_i(p)\}=:J[\varphi_i]$. According to
 Exercise 3G (which applies to arbitrary lattices) we have  $J[\varphi_i]=\sigma(J(L/\theta_i))$, where $\sigma:L/\theta_i\to L$   is the (injective) map that yields the smallest $\varphi_i$-preimages. Hence $\varphi_i$ maps
$J({\cal P}_i)$ bijectively upon $J(L/\theta_i)$.

{\bf 3.4.3} If $J({\cal P}_i)$ is a singleton $\{p\}$ (we also say $p$ is \underbar{projectively isolated}), then $j(L/\theta_i)=|J({\cal P}_i)|=1$ by (10). But then $L/\theta_i=D_2$, and so $p$ is join-prime by Ex.3B. Conversely, if $p\in J(L)$ is join-prime then no  join irreducible $q\neq p$ satisfies $[q_*,q]\approx [p_*,p]$ (Ex.3C), and so
 $\{p\}$ equals some $J({\cal P}_i)$.
In brief: 

\begin{enumerate}
	\item[(11)] {\it In a  modular lattice the projectively isolated join irreducibles are\\ exactly the join-primes.}
\end{enumerate}

 This further implies:
 
\begin{enumerate}
	\item[(12)] {\it In each  distributive lattice $L$ it holds that $s(L)=\delta(L)=j(L)$. In particular $L$ is finite.} 
\end{enumerate}

{\it Proof of (12).}  In a  distributive lattice {\it all} join irreducibles are join-prime (Ex. 2K). Hence (11) implies that $L/\theta_i\simeq D_2$ for all $s(L)$ many co-atoms
$\theta_i\in Con(L)$. From $j(D_2)=1$ and (10) follows $ j(L)=\sum_{i=1}^{s(L)} j(L/\theta_i)=1+\cdots+1=s(L)$. 
Likewise from $\delta(D_2)=1$ and Thm.3.3 follows $\delta(L)=\sum_{i=1}^{s(L)}\delta(L/\theta_i)=1+\cdots+1=s(L)$. Any (FL) lattice $L$ with $j(L)<\infty$ is finite
(of cardinality $\le 2^{j(L)}$) because $x\neq y\Ra J(x)\neq J(y)$. $\square$

A direct proof of $\delta(L)=j(L)$ in (12) is given in Exercise 3A. For non-modular lattices the identity $s(L)=\delta(L)=j(L)$ does not imply that $L$ is distributive (take $L=N_5$).
However, recall from Exercise 2.O  that 'all join irreducibles are join-prime' {\it is} indeed equivalent to distributivity.

\begin{center}
	\includegraphics[scale=0.96]{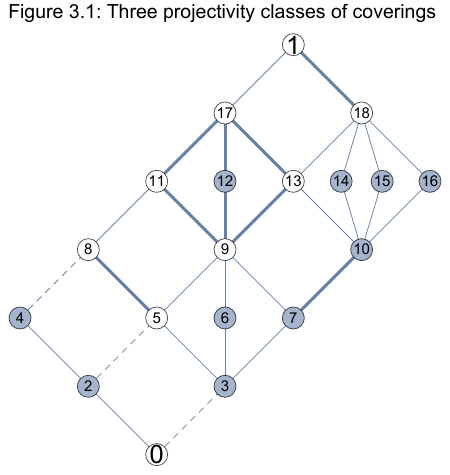} 
\end{center}

{\bf 3.5} Assume the diagram of the modular lattice $L$ is planar (this is not strictly necessary but certainly helps) and that the projectivity classes of coverings are known. Then all factor lattices $L/\theta$ can be drawn by "inspection". 
For instance consider the modular lattice $L_1$ in Figure 3.1. It has three projectivity classes of coverings, which are rendered by {\it n}ormal, {\it d}ashed, and {\it b}oldface lines. Let us denote the sets of join irreducibles belonging to these projectivity classes by $J(n),J(d),$ and $J(b)$ respectively. One checks that

$(13)\qquad J(L_1)=J(n)\uplus J(d)\uplus J(b)=\{2,6,7,14,15,16\}\uplus\{3\}\uplus \{4,10,12\}$

Hence $3$ must be prime by (11). Let $\theta_n,\theta_d,\theta_b$ be the congruences that separate exactly the normal, dashed and boldface coverings respectively.
They are hence the co-atoms of $Con(L_1)$.  To fix ideas let $L=L_1$ and $\theta=\theta_n$. Then the diagram of $L':=L_1/\theta_n$ in Fig. 3.1(B) is obtained by collapsing all dashed and boldface coverings in Fig. 3.1(A). All normal (=normally drawn) coverings thus survive but some may be forced to coincide. For instance (5,9) and (8,11) in Fig. 3.1(A) coincide and become $(p_1,b)$ in Fig.3.1(B). By (10), if $f:L_1\to L'$ is the natural epimorphism then  $f$ maps $J(n)$ bijectively  upon $J(L')$, i.e.

$f(2)=p_1,\ f(6)=p_2,\ f(7)=p_3,\ f(14)=p_4,\ f(15)=p_5,\ f(16)=p_6.$

While most pre-images are singletons, note that $f^{-1}(p_1)=\{2,4,5,8\}=[2,8]$.  Apart from [2,8], the other nontrivial blocks of $\theta_n$ are [0,3], [9,17], and [18,1].

\begin{center}
	\includegraphics[scale=0.7]{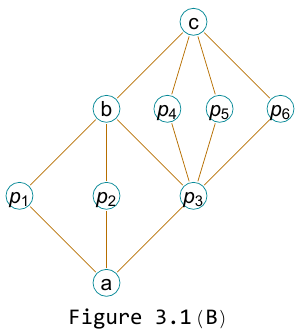} 
\end{center}

As another example, let us look ahead to the modular lattice $L_2$ in Figure 4.1. It also has three projectivity classes of coverings. One checks that the three c-simple factors of $L_2$ are
$M_3, D_2,D_2$. Let us illustrate the statement about maximal 0,1-chains made in Theorem 3.3.
The maximal 0-1-chain $0\prec 2\prec 15\prec 1$ is of type $dbbn$, whereas $0\prec 4\prec 13\prec 19\prec 1$ is of type $bnbd$. 

Do we really need to determine the projectivity classes of coverings in order to find (or just count) the c-simple factors of a  modular lattice? Fortunately there is a faster way. We  handle the atomistic case in Section 4.5 and the general case in Theorem 7.3.

{\bf 3.6 Exercises}

{\bf Exercise 3A.}
\begin{enumerate}
	\item[(i)] Let $(a,b)$ be a covering in any lattice. Show there is $p\in J(a,b)$ with $[p_*,p]\nearrow [a,b]$.\\ (Hint: Consider an element $p\in L$ minimal with the property that $p\le a$ but $p\not\le b$.)
	\item[(ii)] Trimming (1) in Section 2 show that in each  distributive lattice $L$ one has $\delta(L)=j(L)$. (Hint: Use Ex.2A and the join-primeness of all $p\in J(L)$)
	\item[(iii)] Show that  $j(L)\ge s(L)$ in every lattice $L$.
	\item[(iv)] Suppose that the lattice $L$ is such that $(p_*,p)\approx (q_*,q)$ for all $p,q\in J(L)$. Show that $L$ is congruence-simple.
     
\end{enumerate} 

{\it Solution to Exercise 3A.} (i). The set $S\s L$ of all $x$ with $x\le b,\ x\not\le a$ is not empty since $a\in S$. If $x\in S$, then $a<a+x\le b$, and so $a+x=b$. Since $L$ has finite length, $S$ has  minimal members $p$. Suppose $p=u+v$ for some $u,v<p$. By the minimality of $p$ one has $u,v\le a$, yielding the contradiction $p\le a$. It follows that $p\in J(L)$ and $p_*\le a$ and $[p_*,p]\nearrow [a,b]$.

(iii). By (i) each homogeneous projectivity class of coverings contains at least one member of type $[q_*,q]$. Hence $j(L)\ge s(L)$.

(iv). By (i) each $(a,b)\in Cov(L)$ transposes down to some $(p_*,p)$. By assumption any two $(a,b),(c,d)\in Cov(L)$ are therefore projective. Every nonzero congruence $\theta$ contains at least one covering $(a,b)$. According to (2) it follows from $(a,b)\in\theta$ and $(a,b)\approx (c,d)$ that $(c,d)\in\theta$. Hence $Cov(L)\s\theta$. It follows from (4) that $(0_L,1_L)\in\theta$, and so $L$ is congruence-simple. (A harder question is, whether the converse holds. It e.g. holds for modular lattices.)

{\bf Exercise 3B.} Prove the claims about join-primes made in Section 3.1.

{\bf Exercise 3C.} Let ${\cal Q}\s Cov(L)$ be closed under sub-transpositions. Show the hard part of Lemma 3.1, i.e. that $Q=Cov(\theta)$ for some $\theta\in Con(L)$.

{\it Solution to Exercise 3C.} Given ${\cal Q}$ define the binary relation $\theta$ by setting $a\theta b$ iff there is a chain from $ab$ to $a+b$ whose coverings  are from ${\cal Q}$. In order to show that $\theta$ is a congruence it suffices to establish (3),(4),(5), the first two of which being straightforward.
As to (5), let $a<b$ be such that $a\theta b$, and pick $c\in L$ such that $ca<cb$ (for $ca=cb$ the claim is trivial). In order to show\footnote{Showing $(c+a)\theta(c+b)$ proceeds along the same lines and will be omitted.} that $(ca)\theta(cb)$, consider any maximal chain $x_0\prec x_1\prec\cdots\prec x_{n-1}\prec x_n$ such that $x_0=a,\ x_n=b$ and $(x_i,x_{i+1})\in{\cal Q}$ for all $0\le i<n$. Putting $y_0:=cx_0$ we will construct a chain
$y_0<y_1<\cdots <y_s$ such that $y_s=cx_n$ and for each interval $[y_j,y_{j+1}]$ there is a covering  $[x_i,x_{i+1}]$ that transposes down to it. Because by assumption ${\cal Q}$ is closed under sub-transposition, this implies that the coverings  of any maximal chain from $y_j$ to $y_{j+1}$ belong to ${\cal Q}$. Connecting these $s$ maximal chains yields a maximal chain from $ca$ to $cb$ all of whose coverings  belong to ${\cal Q}$, and so $(ca)\theta(cb)$.

For the sake of notation (and because the argument can be repeated) we only explain how $y_1$ and $y_2$ are found. Because of $cx_0<cx_n$ there is an index $i$ such that $cx_0=cx_1=\cdots=cx_i<cx_{i+1}$. Putting $y_1:=cx_{i+1}$ we have $x_i y_1=x_i c=y_0$. Furthermore $y_1\not\le x_i$, hence $x_i<x_i+y_1\le x_{i+1}$, hence
$x_i+y_1=x_{i+1}$ (in view of $x_i\prec x_{i+1}$), hence $[x_i,x_{i+1}]\searrow  [y_0,y_1]$. If $y_1=cx_n$, we are done (put $s:=1$). Otherwise $y_1<cx_n$, and so there is an index $k$ with $cx_{i+1}=cx_{i+2}=\cdots=cx_k<cx_{k+1}$. Putting $y_2:=cx_{k+1}$ one similarly concludes that $[x_k,x_{k+1}]\searrow  [y_1,y_2]$.

{\bf Exercise 3D.} Use Lemma 3.1 to calculate the congruence lattice $Con(N_5)$.

{\bf Exercise 3E.}  In each of these modular lattices, determine the projectivity classes $\cal P$ of coverings, and the corresponding c-simple factor lattices obtained by separating exactly the coverings in $\cal P$.

\begin{center}
	\includegraphics[scale=0.7]{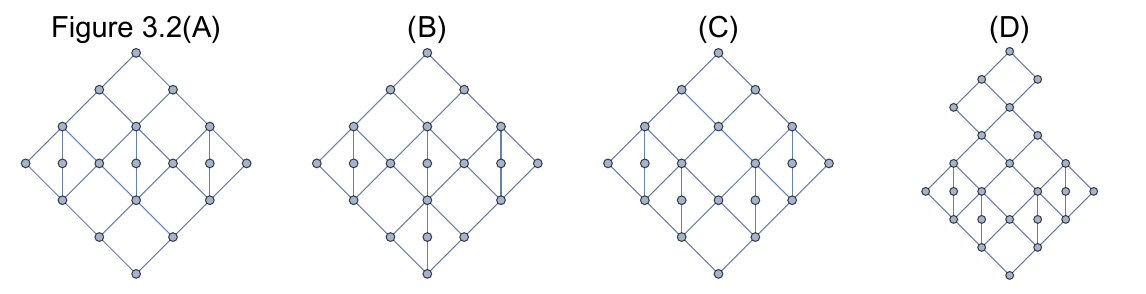} 
\end{center}

{\bf Exercise 3F.} Construct a lattice with two coverings $(a,b),(c,d)$ which are homogeneously projective. Yet  there also is a witnessing sequence for $(a,b)\approx (c,d)$ which features an interval of length two.

{\bf Exercise 3G.} Let $f:L\to L'$ be an epimorphism between any two lattices. Define $\sigma:L'\to L$ by $\sigma(y):=min f^{-1}(y)$ (=smallest pre-image of $y$). Clearly $\sigma$ is injective. Further show:
\begin{enumerate}
	\item[(i)] $\sigma(J(L'))\s J[f]:=\{p\in J(L):\ f(p_*)<f(p)\}$
	\item[(ii)] $\sigma$ is monotone
	\item[(iii)] $\sigma$ is a $\vee$-embedding (see Sec. 2.2)
	\item[(iv)] Improving (i) show that $\sigma(J(L'))= J[f]$.
\end{enumerate}

{\it Solution to Exercise 3G.} (i) Take $p'\in J(L')$ and suppose there were $a,b<\sigma(p')$ with $a+b=\sigma(p')$. Then $f(a)<p',\ f(b)<p'$, and so 
$p'=f(\sigma(p'))=f(a+b)=f(a)+f(b)$. This contradicts  $p'\in J(L')$. Hence $\sigma(J(L'))\s J(L)$. From $\sigma(p')=:p\in J(L)$ follows $f(p_*)<f(p)$, and so $\sigma(J(L'))\s J[f]$.

(ii) Let $y'<z'$ be in $L'$ and suppose we had $\sigma(y')\not<\sigma(z')$. Then $\sigma(y')\wedge\sigma(z')<\sigma(y')$, and so $f(\sigma(y')\wedge\sigma(z'))<y'$. This contradicts $f(\sigma(y')\wedge\sigma(z'))=f(\sigma(y'))\wedge f(\sigma(z'))=y'\wedge z'=y'$.

(iii) Take $x',y'\in L'$. Since $\sigma$ is monotone, we have $\sigma(x')\le\sigma(x'\vee y')$ and $\sigma(y')\le\sigma(x'\vee y')$, which yields $\sigma(x')\vee \sigma(y')\le\sigma(x'\vee y')$. On the other hand $f(\sigma(x')\vee\sigma(y'))=f(\sigma(x'))\vee f(\sigma(y'))=x'\vee y'$ yields $\sigma(x'\vee y')\le\sigma(x')\vee\sigma(y')$. Hence $\sigma(x'\vee y')=\sigma(x')\vee\sigma(y')$.

 (iv) Take any $p\in J[f]$. Putting $x':=f(p)$ we show (a) that $p$ is the smallest $f$-preimage of $x'$, and (b) that $x'\in J(L')$. This  establishes the claim $J[f]\s\sigma(J(L'))$. As to (a), by (1) the full preimage $S$ of $x'$ is an interval of $L$ with $p\in S$. Let $s$ be the smallest element of $S$. From $s<p$ would follow $s\le p_*<p$, hence
$f(s)\le f(p_*)<f(p)$, contradicting $f(s)=f(p)$. Hence $p=s$. As to (b), suppose we had $x'=a'\vee b'$ for some $a',b'< x'$. From (iii) we get $p=\sigma(x')=\sigma(a'\vee b')=\sigma(a')\vee\sigma(b')$, which contradicts $p\in J(L)$. The contradiction shows that $x'\in J(L')$.

\bigskip
\centerline{\bf 4. Geometric lattices and matroids}
\medskip

In a nutshell, we start out with a subclass of semimodular lattices, i.e. geometric lattices. This leads us to matroids, their connected components,  matroid-circuits, and coordinatizing matroids by fields. Much in this Section is inspired by Faigle's chapter in [Wh1].

{\bf 4.1} What are two natural necessary conditions for any   lattice to be sectionally complemented? The first, trivially, is to be {\it complemented}. The second is to be {atomistic}. Theorem 2.2  shows that for modular lattices either one of them is sufficient.  For {\it semi}modular lattices 'complemented' is not sufficient, as evidenced by $SM_{10}$ in Fig. 2.1,
but 'atomistic' does the job:

{\bf Lemma 4.1 }{\it  Let $L$ be semimodular. Then $L$ is sectionally complemented iff it is atomistic.}

{\it Proof.} It remains to show that a semimodular, atomistic lattice $L$ is sectionally complemented. So why is each  interval $[a,b]\s L$ complemented? Since $[a,b]$ is itself semimodular, it suffices by Ex.2L to show that $b$ is a sum of atoms of $[a,b]$. By assumption there are atoms $p_1,...,p_t$  of $L$ with $p_1+\cdots +p_t=b$. Evidently
$(a+p_1)+\cdots +(a+p_t)=b$. Now $a+p_i=a$ if $p_i\le a$. Otherwise it follows from $[0,p_i]\nearrow [a,a+p_i]$
and semimodularity that $a\prec a+p_i$. Hence $b$ is a sum of atoms of $[a,b]$. $\square$

 By definition a lattice is \underbar{geometric} if it is semimodular and satisfies the two equivalent conditions in Lemma 4.1. Each   lattice $L$ is isomorphic to a direct product $L_1\times\cdots\times L_t$ of directly irreducible lattices $L_i$.
It is routine to show (Ex.4A) that $L$ is  geometric iff all factors $L_i$ are geometric.  
According to Ex.2L each geometric lattice is co-atomistic.

{\bf 4.2} The closure space $\M=(E,cl)$ is called a \underbar{matroid} if $\La(\M):=\La(E,cl)$ has finite length and  $cl$ satisfies the \underbar{Exchange Axiom}, i.e. for all $p,q\in E$ and $X\s E$ it holds that
\begin{enumerate}
	\item[(1)]{\it \qquad  $\ p\not\in cl(X)$ and $p\in cl(X\cup\{q\})\ $ jointly imply $q\in cl(X\cup\{p\})$.} 
\end{enumerate}

Condition (1) has a crisper equivalent: Whenever an atom $q$ is not in a closed set $cl(X)$, then $cl(X\cup\{q\})$ is an upper cover of $cl(X)$ in the lattice   $\La(\M)$. In particular, looking at the zero $cl(\es)$ of $\La(\M)$, we see that the sets $[p]:=cl(\{p\})$, where $p$ ranges over $E\setminus cl(\es))$, are exactly\footnote{ The elements in $cl(\es)$ are called \underbar{loops}, and the non-loops in $[p]$ are \underbar{parallel} to $p$. (The terminology comes from graph theory, see 4.2.3.) One calls	$\M$ \underbar{simple} if $cl(\es)=\es$ and $[p]=\{p\}$ for all $p\in E$. This kind of 'simple'  has nothing to do with c-simple lattices.} the atoms of $\La(\M)$.
  It readily follows (Ex. 4B) that $\La(\M)$ is a geometric lattice. The members of $\La(\M)$ are called the \underbar{flats} of $\M$.
  
 {\bf 4.2.1.} An important type of matroid is obtained as follows. For a field $k$ and $n\in \N$  fix a (possibly infinite) subset  $E'\s \La(k^n)$ of vectors. Since for $X\s E'$ the generated subspace $\langle X\rangle$ need not be contained in $E'$, define $cl'(X):=\langle X\rangle\cap E'$. In this scenario (1) boils down to the Steinitz Exchange Lemma familiar from  linear algebra, and so $(E',cl')$ is a matroid, called a \underbar{$k$-linear} matroid. It is usually hard to decide, whether some given matroid $(E,cl)$ is \underbar{$k$-coordinatizable}, i.e.  isomorphic\footnote{Generally an \underbar{isomorphism} $f:E\to E'$ between matroids $(E,cl)$ and $(E',cl')$ by definition satisfies $f(cl(X))=cl'(f(X))$ for all $X\s E$.}  to a suitable $k$-linear matroid $(E',cl')$.

{\bf 4.2.2.} A subset $I\s E$ is \underbar{independent} in a matroid $\M=(E,cl)$ if $x\not\in cl(I\setminus\{x\})$ for all $x\in I$. Each subset of an independent set is independent. Furthermore (Ex. 4C): 
 
 \begin{enumerate}
 	\item[(2)] {\it If $I\s E$ is independent, 
 		and $p\not\in cl(I)$, then $I\cup\{p\}$ stays independent.  } 
 \end{enumerate}
 
Let $X$ be a flat of $\M$ and let $I\s X$ be a maximal independent subset. Such $I$ cannot be infinite (Ex.4D), and by (2) there is no $p\in X\setminus cl(I)$. Hence $cl(I)=X$. We see that  all maximal independent subsets of a flat generate the flat. Further  they have the {\it same} finite cardinality (which coincides with the height $\delta(X)$ of $X$ in $\La(\M)$).
More generally, the \underbar{rank} of any subset $Y\s E$ is defined as $rk(Y):=\delta(cl(Y))$.
 A maximal independent set $B\s E$ is called a \underbar{base} of $\M$. Restating the above for  $X=E$ we record:
 
  \begin{enumerate}
 	\item[(3)] {\it All bases of a matroid $(E,cl)$ have cardinality $rk(E)$.  } 
 \end{enumerate}

{\bf 4.2.3} Each minimal \underbar{dependent} (=not independent) set $C$ is called a \underbar{circuit} of $\M$. A priori this merely implies there is at least {\it one} $x\in C$ with $x\in cl(C\setminus\{x\})$. However, it follows from (2) that $x\in cl(C\setminus\{x\})$ for {\it all} $x\in C$.
 One can retrieve $cl$ from the set $\Ci$ of all circuits via the implicational base\footnote{A set of type $C\setminus\{x\}$ is called a \underbar{broken circuit}. Thus $X\s E$ is a flat iff with any broken circuit it contains the whole circuit.}
 
 $(4) \qquad \Sigma:=\{\ (C\setminus \{x\})\to\{x\}:\ C\in \Ci,\  x\in C\}.$
 
 Not only can one retrieve $cl$ from the circuits, one can also start out with the circuits in the first place.  Specifically, let $\Ci$ be any family of subsets (called 'circuits') of a set $E$, such that two axioms are satisfied.
 
 \begin{itemize}
 	\item {\it Axiom 1:}  $\Ci$ is an antichain.
 	\item {\it Axiom 2:} For all intersecting  $C,C'\in\Ci$ and all $x\in C\cap C'$\\ there is a $C''\in\Ci$ with $C''\s (C\cup C')\setminus\{x\}$.
 	
 \end{itemize}

   Then it holds that the closure space $(E,cl)$ defined by the implications in (4) is a matroid. We will also write $\M(E,\Ci)$.
 For instance if $G=(V,E)$ is a graph with vertex set $V$ and edge set $E$, then the circuits of $G$ (in the usual sense) constitute a set family $\Ci$ that satisfies  Axiom 1 and Axiom 2. One calls $\M(E,\Ci)$ the \underbar{graphic matroid} induced by $G$.

{\bf 4.2.4.} Let $\M=(E,cl)$ be a matroid. Two elements of $E$ are $\M$-\underbar{connected} if there is a circuit containing both of them. From Axiom 2  follows that 'being connected' is an equivalence relation. Its blocks, i.e. the \underbar{connected components} $A_1,\ldots, A_s$, are closed (so $B\s A_i\Ra cl(B)\s A_i$) and have the property that each $X\s E$ can be handled componentwise:

$(5)\qquad cl(X)=cl(X\cap A_1)\uplus\cdots\uplus cl(X\cap A_s).$

Attention, if $G$ is a graph, then the $\M(G)$-connected sets are the edge sets of the 2-{\it connected} components of $G$. Yet another concept of 'connected' will be introduced in Section 4.4.

{\bf 4.2.5} The direction  can be reversed. Namely, putting $s=2$ for convenience, let $\M_1,\M_2$ be matroids on {\it disjoint} sets $A_1, A_2$, having closure operators $cl_1,cl_2$.
Putting $A:=A_1\cup A_2$ one gets a new matroid $\M=(A,cl)$ with $cl$ defined by

$\quad cl(X):=c_1(A_1\cap X)\cup c_2(A_2\cap X).$

One calls $\M$ the \underbar{direct sum} of $\M_1$ and $\M_2$.

{\bf 4.3} Recall that $\La(\M)$ is geometric for each matroid $\M$. Here comes  a kind of converse. For any geometric lattice $L$ define the closure operator $cl_L$ on the set $At(L)=J(L)$ of atoms  by

 $(6)\qquad cl_L(X):=J(\sum X)$ {\it for all} $X\s J(L)$.

It follows (Ex. 4E) that
 $\M(L):=(At(L),cl_L)$ a simple matroid and that $\La(\M(L))\simeq L$. We prefer to work with  simple matroids of type $\M(L)$ because then points and atoms are the same thing.

{\bf 4.3.1} Join irreducibles $p$ and $q$ in a  lattice are called \underbar{perspective} if the coverings  $(p_*,p)$ and $(q_*,q)$ have a common upper transpose. For instance, the join irreducibles $5$ and $9$ in $SM_{10}$ are perspective: $(2,5)\nearrow (7,1)\nwarrow (4,9)$. In semimodular lattices the common upper transpose must be a covering. Therefore perspectivity is most useful in semimodular lattices, particularly  geometric or modular ones. 

{\bf Lemma 4.2} {\it Let $L$ be a geometric lattice. 
\begin{enumerate}
	\item[(a)] 	Any  $p\neq q$ in $At(L)$ are perspective iff they belong to the same connected component of $\M(L)$.
\item[(b)] If additonally $L$ is modular, then $p\neq q$ are perspective iff there a third atom $r$ 
such that $\{0,p,q,r,p+q\}$ is a covering $M_3$-sublattice. 

\end{enumerate}	}

{\it Proof.} (a) By 4.2.4 we need to show that being in a common circuit of $\M(L)$ is equivalent to being perspective. First let $p,q$ be distinct elements of any circuit $C\s At(L)$. Let $e$ be the sum of all atoms in $C\setminus \{p,q\}$. We conclude:

$(7)\quad p\le e+q\ {(1)\atop \Ra}\ q\le e+p\ \Ra\ e+p=e+q=:f\ \Ra\ [0,p]\nearrow [e,f]\searrow [0,q].$

Conversely, if the atoms $p\neq q$ are perspective, we choose the common upper transpose $[e,f]$ of $[0,p]$ and $[0,q]$ such that $\delta(e)$ is minimal. By 4.2.2 there are independent atoms $r_1,\ldots, r_t$  with $r_1+\cdots r_t=e$. Put $C:=\{r_1,\ldots, r_t,p,q\}$. From $q<f=e+p=r_1+\cdots r_t+p$ follows the dependency of $C$. It remains to show that dropping any element from $C$ yields an independent set.
From $p\not\le r_1+\cdots r_t\ (=e)$ and (2) follows that $\{r_1,\ldots, r_t,p\}$ is independent. Likewise $\{r_1,\ldots, r_t,q\}$ is independent. It remains to check, say, $\{r_2,\ldots, r_t,p,q\}$.
Put $e':=r_2+\cdots + r_t$. If we had $p\le e'+q$, then also $q\le e'+p$ by (1). Putting $f':=e'+p=e'+q$,
this yields the contradiction $[0,p]\nearrow [e',f']\searrow [0,q]$ with $\delta(e')<\delta(e)$. 

(b) By (a) and Ex.4F, in any geometric lattice atoms $p,q$ are perspective iff they have a common complement, i.e.
there is $c\in L$ 
with $p+c=q+c=1$ and $pc=qc=0$. Putting $r:=c(p+q)$ we conclude (using the modularity assumption) that
$p+r=p+c(p+q)=(p+c)(p+q)=p+q$ and $pr=pc(p+q)=0$. By symmetry also $q+r=p+q,\ qr=0$. It is clear that $r$ is an atom distinct from $p,q$ and that  $\{0,p,q,r,p+q\}$ is a covering $M_3$-sublattice.
Thus in the modular case the common upper transpose $(e,f)$ in (7) can be found as close to $(p_*,p),(q_*,q)$ as possible, i.e $(e,f)=(r,p+q)$.
$\square$

{\bf 4.4} In any  lattice perspective join irreducibles $p,q$ are a fortiori \underbar{projective} (written $p\approx q$) in the sense (Sec.3.4.2) that $(p_*,p)\approx (q_*,q)$. The converse fails (Ex.4G), but holds in geometric lattices.

\begin{enumerate}
	\item[(8)\quad] {\it In any geometric lattice $L$ any two atoms are projective iff they are perspective. Further, let $A_1,...,A_s\s At(L)$ be the projectivity classes of atoms, and $a_i:=\sum A_i$. Then $L\simeq [0,a_1]\times\cdots\times [0,a_s]$. All lattices $[0,a_i]$ are c-simple.}
	\end{enumerate}

The direction from right to left in the claimed isomorphism is given by $(x_1,...,x_s)\mapsto x_1+\cdots +x_s$, whereas the converse direction is
$x\mapsto (x_1,...,x_s)$ where $x_i:=\sum\{p\in J(x):\ p\le a_i\}$. While in every lattice "c-simple $\Ra$ directly irreducible", by (8) the two are equivalent for geometric lattices.

{\it Proof of (8)}. The first claim is referred to Exercice 4H. If $A_1,...,A_s\s E$ are the connected components of any matroid  $\M=(E,cl)$, then $\La(\M)\simeq
\La(\M(A_1))\times\cdots\times \La(\M(A_s))$ by Ex.4I. If $E$ happens to be the set of atoms of $L$, this isomorphism becomes $L\simeq [0,a_1]\times\cdots\times [0,a_s]$. It follows from Ex.3A(iv) that each $[0,a_i]$ is c-simple. $\square$

{\bf 4.5} At the end of Section 3.5 we raised the question whether finding, or at least counting, the c-simple factors of a  modular lattice can be achieved without determining all projectivity classes ${\cal P}_i$  of coverings, and thus the projectivity classes $J({\cal P}_i)$ of join irreducibles. Specializing to {\it geometric} modular lattices we can now provide a smooth solution. (In view of (8) we can also say that we specialize to {\it modular} geometric lattices.) Thus by (8) 'projective' boils down to 'perspective', which by Lemma 4.2(b) amounts to  easily detect or exclude a $M_3$ sublattice. For instance take $L$ from Figure 4.1.
 The atoms $3$ and $4$ are 
perspective since  $[0,3+4]\simeq M_3$, whereas $2$ and $3$ are not perspective since $[0,2+3]\not\simeq M_3$. It is clear that the perspectivity (=projectivity) classes of atoms are $\{2\},\{6\}\{3,4,5\}$,
and so by (8) $L$ has three c-simple factors $[0,a_1]\simeq D_2,\ [0,a_2]\simeq D_2,\ [0,a_3]\simeq M_3$, where $a_1=2,\ a_2=6,\ a_3=3+4+5=10$.

\begin{center}
	\includegraphics[scale=0.9]{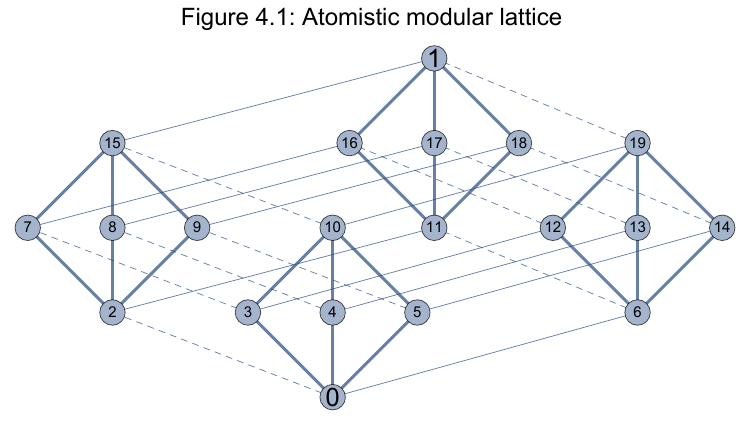} 
\end{center}

In fact, if $L$ in (8) is modular, then the lattices $[0,a_i]$ are isomorphic to subspace lattices of so-called projective spaces; more on that in 5.6.

{\bf 4.6 Exercises}

{\bf Exercise 4A.} Consider the direct product of  lattices $L=L_1\times L_2$. Show that $L$ is geometric iff $L_1$ and $L_2$ are geometric.

{\bf Exercise 4B.} Show that each matroid $\M=(E,cl)$ has a geometric flat lattice $\La(\M)$.

{\it Solution of Ex. 4B.} For all $X\in L(M)$ one has $X=cl(\bigcup \{\ [p]:\ p\in X\})=\bigvee\{[p]: p\in X\}$, and so $L(M)$ is atomistic. To establish semimodularity, take any $A,B,C,D\in L(M)$ with $A\prec B$ and $[A,B]\nearrow [C,D]$. We need to show $C\prec D$. For any $p\in B\setminus A$ we have $cl(A\cup\{p\})=B$ because of $A\prec B$. Hence $D=B\vee C=cl(B\cup C)=cl(A\cup\{p\}\cup C)=cl(C\cup\{p\})$. As noted after (1), this implies $C\prec D$.

{\bf Exercise 4C.} Prove property (2). Then use (2) to show that for each circuit $C$ it holds that $x\not\in C\setminus\{x\}$ for {\it all} $x\in C$.

{\it Solution of Ex. 4C} By way of contradiction, suppose $X\s E$ was independent, $p\not\in cl(X)$, but $X\cup\{p\}$ was dependent. Hence some element in this set (by assumption not $p$) is in the closure of the others: $q\in cl((X\setminus \{q\})\cup \{p\})$. But  $q\not\in cl(X\setminus \{q\})$ since $X$ is independent.
From (1) follows the desired contradiction  $p\in cl((X\setminus \{q\})\cup \{q\})=cl(X)$. 

{\bf Exercise 4D.} Show that  a matroid cannot contain infinite independent sets.

{\it Solution of Ex. 4D} Suppose $\M=(E,cl)$ had a (w.l.o.g. denumerably) infinite independent set $\{p_1,p_2,\ldots\}\s E$. By definition of 'independent' $p_{i+1}\not\in cl(\{p_1,..,p_i\})$, and so we had an infinite chain of flats $cl(\{p_1\})\subset cl(\{p_1,p_2\}))\subset\cdots$, contradicting the fact that $\La(\M)$ has finite length by definition of 'matroid'.

{\bf Exercise 4E.} Show that the closure space $\M(L)=(At(L),cl_L)$ associated to the geometric lattice $L$ is a simple matroid.

{\it Solution of Ex. 4E}
In order to establish the exchange property (1), let $X\s At(L)$ and let $p\neq q$  be atoms not in $cl(X)=:x$, but satisfying $p\in cl(X\cup\{q\})$. Hence $p\le x+q$, and so $x<x+p\le x+q$. But $x\prec x+q$ by semimodularity. Hence $x+p=x+q$, which implies $q\in cl(X\cup\{p\})$.

{\bf Exercise 4F.}  Show that in any  lattice $L$ two atoms $p,q$ which possess a common complement are perspective. If $L$ is co-atomistic (e.g. geometric), the converse holds as well.

{\it Solution of Ex.4F.} If $c$ is a common complement, then $[0,p]\nearrow [c,1]\nwarrow [0,q]$. Conversely, let $[0,p]\nearrow [e,f]\nwarrow [0,q]$.
If $L$ is co-atomistic, there is a co-atom $c$ with $e\le c$ but $f\not\le c$. In view of $p+e=f$ this implies $p\not\le c$. Hence $c$ is a complement of $p$. Likewise $c$ is a complement of $q$.

{\bf Exercise 4G.} Find two  join irreducibles in $L_1$ of Fig. 3.1 which are projective but not perspective.

{\it Solution of Ex.4G.} For instance $7,16\in J(L_1)$ are projective, witnessed by $(3,7)\nearrow(6,9)\searrow (0,2)\nearrow (14,18)\searrow (10,16)$. The atoms $7,16$ are not perspective since generally coverings $(a,b),(c,d)$ with $b\le c$ cannot have a common upper transpose (check).

{\bf Exercise 4H.} Show that in a geometric lattice two atoms are perspective iff they are projective.

{\it Solution of Ex.4H.} The reverse implication being trivial we only show 'projective $\Ra$ perspective'. So fix $p,q\in At(L)$ with $p\approx q$. Since neither $(0,p)$ nor $(0,q)$ can be upper transposes, each witnessing sequence of $p\approx q$ is of type

$(0,p)\nearrow (a,b)\searrow (c,d)\nearrow\cdots\searrow (e,f)\nearrow (g,h)\searrow (0,q).$

By Ex.3A(i) each lower transpose $(c,d),...,(e,f)$ can be assumed to stem from an atom (so $(c,d)=(0,r)$, etc).
Therefore $p\sim r\sim\cdots\sim q$. Since $\sim$ is transitive by Lemma 4.2(a), we conclude that $p\sim q$.

{\bf Exercise 4I.} Let $A_1,...,A_s$ be the connected components of $\M=(E,cl)$. Using (5) show that
$\La(\M)\simeq\La(\M(A_1))\times\cdots\times \La(\M(A_s))$.

{\bf 5. Partial linear spaces and matroids}

In 5.1 and 5.2 we  put in place standard material about partial linear spaces (PLS), a general reference being [BB]. Some technicalities about cycles (5.3) seem to be new. Also novel, but perhaps more exhilarating,  is 
the concept of a 'point-splitting' (5.4).
In  5.6 (projective geometries) and 5.7 (parallel connections of lines) PLSes are related to matroids.

{\bf 5.1} Let $E$ be a set, whose elements we call \underbar{points}, and let $\Lambda\s \Pow(E)$ be a set of subsets called \underbar{lines}. The pair  $(E,\Lambda)$ is a {\it partial linear space (PLS)} if the following is satisfied [BB].

\begin{enumerate}
		\item[\qquad (PLS 1)]{\it Every line is incident with at least two points.}
		\item[\quad (PLS 2)] {\it Any two distinct points are incident with at most one line.}
		\end{enumerate}

  Notice that axiom (PLS 2) is equivalent to
 
 \begin{enumerate}
\item[(PLS 2')] For all distinct $\ell,\ell'\in \LL$ it holds that $|\ell\cap \ell'|\le 1$. 
\end{enumerate}

When distinct points $p,q\in E$ are incident with a line, then it is unique by (PLS 2) and we denote it by $\langle p,q\rangle$. A line $\ell\in\LL$  in a PLS $(E,\LL)$ is \underbar{short} if $|\ell|=2$, otherwise it is \underbar{long}. If all lines in a PLS are short, then the PLS is just a graph. In the PLS below (Fig.5.1(A)), the so-called \underbar{Fano plane}, all seven lines qualify as long, albeit "just quite". (It is impossible to display all 7 lines as 'straight lines'.)

{\sl Figure 5.1: Making the Fano plane acyclic}

\begin{center}	\includegraphics[scale=0.82, angle=0]{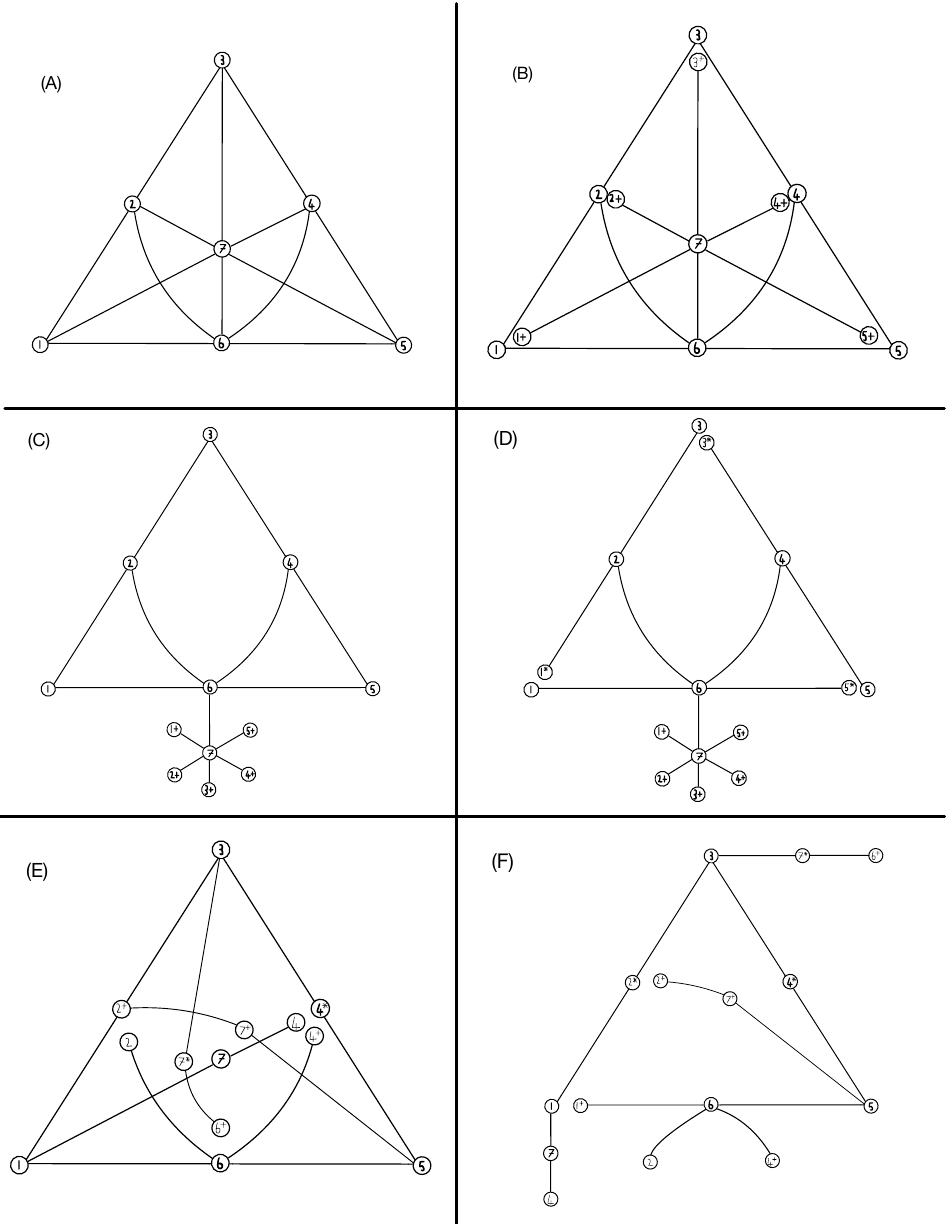} \end{center}

{\sl Figure 5.1: Making the Fano plane acyclic}

{\bf 5.2}
 A \underbar{path} of a PLS $(E,\LL)$ is a sequence $(p_1,p_2,...,p_n)$ of points such that all lines $\langle p_i,p_{i+1}\rangle$ exist $(1\le i<n)$.
Two points $p,q\in E$ are \underbar{connected} if $p=q$ or if they lie on a common path. Let  $E_1,...,E_c$ be the classes of this equivalence relation. Put $\LL_i:=\{\ell\in \LL:\ \ell\s E_i\}$. 
Then $\LL=\biguplus\{\LL_i:\ 1\le i\le c\}$, and $\LL_i=\emptyset$ iff
 $E_i=\{p\}$ consists of an \underbar{isolated} point. The $c=c(E,\LL)$ many PLSes $(E_i,\LL_i)$ are
the \underbar{connected components} of $(E,\LL)$. The PLS $(E,\LL)$ is \underbar{connected}
if $c(E,\LL)=1$.

Apart from connected matroids and connected PLSes let us introduce  yet another meaning of 'connected'. We call a {\it line-set} $S\s \LL$   \underbar{connected} if either $S=\emptyset$ or the PLS $(\bigcup S,S)$ is connected. Note the following:

\begin{enumerate}
	\item[(1)] Let $(E,\LL)$ be a PLS and $S\s\LL$. If $(E,S)$ is a connected PLS, then $S$ is a connected line-set. Conversely, if $S$ is a connected line-set, then $(E,S)$ need not be a connected PLS because it has isolated points whenever $\bigcup S\neq E$.
\end{enumerate}

The next observation is slightly more subtle (Ex.5A):
\begin{enumerate}
	\item[(2)] Each nonempty connected line-set $S$   contains at least one   $\ell\in S$ such that $S\setminus\{\ell\}$ remains a connected line-set.
	\end{enumerate}

{\bf 5.3} A \underbar{cycle} of a PLS $(E,\LL)$ is a sequence $C=(\ell_1,\ell_2,...,\ell_n)$ of $n\ge 3$ distinct lines such that  all\footnote{Without further mention, this means 'modulo $n$', i.e. $1\le i\le n-1$, as well as $\ell_n\cap\ell_1\neq\emptyset$. The distinctness of the  \underbar{junctions} $p_i\in\ell_i\cap\ell_{i+1}$ is crucial. Otherwise unpleasant sequences of lines, such as $(\{1,2,3\},\{1,4\},\{1,5,6,7\})$, would qualify as cycles.} $n$ intersections
$\ell_i\cap\ell_{i+1}$
are nonempty and distinct.
For instance, if the PLS is the Fano plane and\footnote{The  order in which the elements of a line are listed is irrelevant, but some orderings are more illuminating than others. } $\ell_1=\{1,3,2\},\ \ell_2=\{2,7,5\},\ \ell_3=\{5,3,4\},\ 
\ell_4=\{4,7,1\}$, then $(\ell_1,\ell_2,\ell_3,\ell_4)$ is a cycle. 

We define the \underbar{$C$-interior} of a line $\ell_i$ in a cycle $C$ as $\ell_i^\circ:=\ell_i\setminus(\ell_{i-1}\cup\ell_{i+1})$. Clearly consecutive interiors $\ell_i^\circ,\ \ell_{i+1}^\circ$ are disjoint. The cycle $C$ is \underbar{strong} if {\it all} sets $\ell_i^\circ$ are disjoint. For the cycle $C$ above we have $\ell_1^\circ=\ell_3^\circ=\{3\},\ell_2^\circ=\ell_4^\circ=\{7\}$, and so $C$ is not strong.

\begin{enumerate}
	\item [(3)] Each cycle $(\ell_1,...,\ell_n)$ \underbar{contains} a strong cycle in the sense that some subset of $\{\ell_1,...,\ell_n\}$, suitably ordered, constitutes a strong cycle. 
\end{enumerate}

We establish (3) by inducting on $n$. Letting $n=3$ consider a cycle $(\ell_1,\ell_2,\ell_3)$. Let $\ell_1\cap\ell_2=\{p\},\ \ell_2\cap\ell_3=\{q\},\ \ell_3\cap\ell_1=\{r\}$. Hence $\ell_1^\circ=\ell_1\setminus\{p,r\},\ \ell_2^\circ=\ell_2\setminus\{p,q\}$. Assuming $\ell_1^\circ\cap\ell_2^\circ \neq\es$, say $t\in \ell_1^\circ\cap\ell_2^\circ$, yields the contradiction $|\ell_1\cap\ell_2|\ge|\{p,t\}|= 2$. Similarly one shows $\ell_1^\circ\cap\ell_3^\circ=\ell_2^\circ\cap\ell_3^\circ =\es$, and so $(\ell_1,\ell_2,\ell_3)$ is strong. Let $n>3$ and assume that $(\ell_1,...,\ell_n)$ is not strong. Then, wlog, there is $i\in\{3,..,n-1\}$ with $\ell_1^\circ\cap\ell_i^\circ\neq\es$. Then $(\ell_1,\ell_2,..,\ell_i)$ is a cycle of length $<n$, which by induction contains a strong cycle. $\square$

For instance, if the non-strong cycle $(\ell_1,\ell_2,\ell_3,\ell_4)$ above is replaced by $(\ell_1,\ell_2,\ell_3)$, the latter is strong with
$\ell_1^\circ=\{1\},\ell_2^\circ=\{7\},\ell_3^\circ=\{4\}$.

{\bf 5.3.1} If $(E,\LL)$ is a PLS then a line-set $T\s\LL$ is \underbar{acyclic} if no subset can be ordered to become a cycle.

\begin{enumerate}
	\item [(4)] If $C=(\ell_1,...,\ell_n)$ is a strong cycle in some PLS, then $|\ell_1\cup\cdots\cup\ell_n|=|\ell_1|+\cdots+|\ell_n|-n$. \\
	If $S=\{\ell_1',...,\ell_m'\}$ is acyclic and $(\bigcup S,S)$ has $c$ connected components, then
	$|\ell_1'\cup\cdots\cup\ell_m'|=|\ell_1'|+\cdots+|\ell_m'|-m+c$. In particular $|\bigcup S|\ge 2m-m+c=m+c$, and if all $\ell_i'$ are long then $|\bigcup S|\ge 2m+c$. 
\end{enumerate}

As to the first claim, if $Jc$ is the $n$-set of all junctions of $C$, then $|\ell_1\cup\cdots\cup\ell_n|=|\ell_1^\circ\uplus\cdots\uplus\ell_n^\circ\uplus Jc|=(|\ell_1|-2)+\cdots+
(|\ell_n|-2)+n=|\ell_1|+\cdots+|\ell_n|-2n+n$.

As to the second claim, first assume that $c=1$. Then by (2) there is a line in $S$ whose removal does not destroy connectivity. Hence wlog $S_0:=S\setminus\{\ell_m'\}$ is a connected line-set. With $S$ a fortiori $S_0$ is acyclic. Since $S$ is a connected line-set, we have $|(\bigcup S_0)\cap\ell_m'|\ge 1$, But if $\ge 1$ was $\ge 2$ then, since $S_0$ is a connected line-set, $S$ would contain a cycle (which it doesn't). Therefore $\ge 1$ is $=1$. Put another way, $|\ell_m'\setminus (\bigcup S_0)|=|\ell_m'|-1$.
By induction
$|\ell_1'\cup\cdots\cup\ell_{m-1}'|=|\ell_1'|+\cdots+|\ell_{m-1}'|-(m-1)+c$, and so
$|\ell_1'\cup\cdots\cup\ell_m'|=|\ell_1'\cup\cdots\cup\ell_{m-1}'|+(|\ell_m'|-1)=|\ell_1'|+\cdots+|\ell_m'|-m+c$.

Suppose now $c>1$, w.l.o.g. $c=2$. If $(\bigcup S,S)$ has connected components $(S_1,\{\ell_1',...,\ell_m'\})$ and
$(S_2,\{\ell_{m+1}',...,\ell_{m+t}'\})$, then $|S|=|S_1|+|S_2|=(|\ell_1'|+\cdots+|\ell_m'|-m+1)
+(|\ell_{m+1}'|+\cdots+|\ell_{m+t}'|-t+1)=|\ell_1'|+\cdots+|\ell_{m+t}|-(m+t)+c$.

{\bf 5.4} If a PLS $(E,\LL)$ is not acyclic, we can make it acyclic by 'splitting points'.  
Formally a \underbar{point-splitting} is an ordered pair $(\ell,p)\in \LL\times E$ such that $p\in\ell$. Let $PtSp\s \LL\times E$ be the set of all point-splittings. Applying a set $X\s PtSp$ to $(E,\LL)$ yields a new PLS $(E,\LL)^X$, which arises from $(E,\LL)$ as follows. For each $(\ell,p)\in X$ we {\it separate  $\ell$ from $p$}. Let us gradually make this more precise. Let $X':=\{p^\ell:\ (\ell,p)\in X\}$ be  a disjoint $|X|$-element set of 'new points' to be added to $E$. As opposed to the points, the  number of lines $\ell'$ in the new PLS will be the same as the number of old lines $\ell$, due to some natural bijection $\ell\mapsto \ell'$.
Namely, put
$\ell':=(\ell\setminus \{p\in\ell:\ (\ell,p)\in X\})\uplus\{p^\ell:\ (\ell,p)\in X\}.$
The formal definition of the new PLS is

  $(5)\quad (E,\LL)^X:=(E',\LL'),\ where\ E':=E\uplus X'\  and\ \LL':=\{\ell':\ \ell\in\LL\}.$ 

Let us illustrate this on the Fano plane but use $p^*,p^+$ (rather than $p^{\ell_i},p^{\ell_j}$) for the new points resulting from splitting $p$. Referring to Fig. 5.1 (A) we begin by separating the line $\ell :=\{1,7,4\}$ from point $1$. This creates a new point $1^+$ and replaces $\ell$ by $\ell^+:=\{1^+,7,4\}
$. The old point $1$  still belongs to two other lines. We go on and separate $\ell^+$ from $4$. This results in a new point $4^+$ and replaces $\ell^+$ by $\ell^{++}:=\{1^+,7,4^+\}$. Similarly the line $\{2,7,5\}$ in Fig. 5.1(A) gives rise to $\{2^+,7,5^+\}$ in Fig. 5.1(B), and $\{3,7,6\}$ becomes $\{3^+,7,6\}$. Figure 5.1(C) is  a better comprehensible redrawing of Fig. 5.1(B). Continuing with the point-splittings, the lines $\{1,2,3\},\ \{3,4,5\}$ and $\{5,6,1\}$ in Fig. 5.1(C) yield
$\{1^*,2,3\},\ \{3^*,4,5\}$ and $\{5^*,6,1\}$ in Fig. 5.1(D). Evidently the PLS in Fig. 5.1(D) is acyclic and connected.

There are in fact many ways to make the Fano plane acyclic while retaining connectedness. A second one is rendered in Figures 5.1 (E) and (F). Observe that the PLSes in Figures 5.1 (D) and (F) both required 8 point-splittings.

{\bf 5.4.1}
Generally we call $A\s PtSp$ an \underbar{acyclifier} of $(E,\LL)$ if $(E,\LL)^A$ is acyclic but retains the same number of connected components. 

{\bf Theorem 5.1: }{\it Let $\B=(E,\LL)$ be a partial linear space. 
	\begin{enumerate}
		\item [(a)]
	Then all acyclifiers have the same cardinality, which we can hence denote by $r^*(\B)$. 
	\item[(b)] Actually, if $\LL=\{\ell_1,...,\ell_t\}$, then $r^*(\B)=c(\B)+\sum_{i=1}^t|\ell_i|-|E|-t$.
\end{enumerate}
}

For instance, if $\B$ is the Fano plane, then $r^*(\B)=1+7\cdot 3-7-7=8$.

{\it Proof of Thm. 5.1.}  (a) According to [W3,Lemma 11] the acyclifiers of $\B$ bijectively match the bases of
 some\footnote{Each matroid has a \underbar{dual} matroid which is most conveniently defined by its bases; they are the complements of the original bases. Thus above 'some matroid' means 'the dual of some graphic matroid coupled to the PLS'.} matroid, hence by (3) in Sec. 4 all acyclifiers have the same cardinality. 

(b) Once we manage the  connected case (claim (6)), the general case will be immediate.

$(6)\ Each\ connected\ PLS\ \B_1=(E_1,\{\ell_1,...,\ell_m\})\ has\ r^*(\B_1)=1+\sum_{i=1}^m|\ell_i|-|E_1|-m.$

We verify (6) by inducting on $m$. If $m=0$ (so $\B_1$ is an isolated point), then $r^*(\B_1)=0=1+0-1-0$. 
If $m=1$ (so $\B_1$ is a single line), then $r^*(\B_1)=0=1+|\ell_1|-|E_1|-1$. Let $m>1$. Since $\B_1$ is connected and because of (2) there is a line, say $\ell_m$, whose removal results in the connected line-set
$\LL_0:=\LL_1\setminus\{\ell_m\}$. Suppose there are $e\ge 0$ points on $\ell_m$ which are not on any line $\ell_i\neq \ell_m$. Since $\B_1$ is connected, we have $e<|\ell_m|$, and so there are exactly $|\ell_m|-e>0$ points on $\ell_m$ which additionally lie on lines $\neq\ell_m$. By (a) we are free to choose whatever acyclifier suits us best to derive a formula for $r^*(\B_1)$. We opt to launch our point-splitting endeavour by choosing any $|\ell_m|-e-1$ points among the $|\ell_m|-e$ above-mentioned points. By separating $\ell_m$ from these points  we retain connectivity of the PLS and   never again need to split points on $\ell_m$. More formally, let $E_0$ be the subset of $E_1$ obtained by dropping the above-mentioned $e$ points of $\ell_m$. Then $\B_0:=(E_0,\LL_0)$ is a connected\footnote{This is akin to (1): Although $\LL_0$ is a connected line-set, the PLS $(E_1,\LL_0)$ has $e$ isolated points. They disappear when $E_1$ gets replaced by $E_0$.} PLS with fewer lines than $\B_1$. Therefore induction yields $r^*(\B_0)=1+\sum_{i=1}^{m-1}|\ell_i|-|E_0|-(m-1)$. Straightforward arithmetic finishes the proof of (6):

$r^*(\B_1)=(|\ell_m|-e-1)+r^*(\B_0)=(|\ell_m|-e-1)+ 1+\sum_{i=1}^{m-1}|\ell_i|-|E_0|-(m-1)
\\ \hspace*{1.1cm}=1+\sum_{i=1}^{m}|\ell_i|-|E_1|-m$. 

It remains to handle a disconnected $\B$ with $c=c(\B)\ge 2$ connected components $\B_k=(E_k,\LL_k)$. In this case evidently
$r^*(\B)=\sum_{k=1}^c r^*(\B_k)$. From $|E_k|<|E|$ and induction follows that each $\B_k$ is evaluated as $\B_1$ in (6). Summing up these $c$ equations proves claim (b) in  general. $\square$

{\bf 5.5.} A subset  $X\s E$  of a partial linear space $(E,\LL)$ is \underbar{$\LL$-closed}   if

\begin{enumerate}
	\item[(7)\quad]{\it $|X\cap l|\ge 2$ implies $l\s X$, for all $l\in\Lambda$.}
\end{enumerate}

The set $\La(E,\Lambda)$ of all $\LL$-closed subsets is a closure system; viewed as a lattice it is atomistic with
atoms $\{p\}\ (p\in E)$. Let $c_{\LL}:\Pow(E)\to \Pow(E)$ be the coupled closure operator.   A natural
implicational base of  $\La(E,\Lambda)$ (or of $c_{\LL}$) is given by

$(8)\qquad \Sigma_{\LL}:=\{\ \{p,q\}\to \langle p,q\rangle:\ p\neq q\ in\ E\}.$

 Although $\La(E,\Lambda)$  is atomistic, it is generally not semimodular, i.e. not the flat lattice of a matroid. However, in 5.6 (projective spaces) and 5.7 (parallel connection of lines) we construct types of PLSes where it works.

 {\bf 5.6} A partial linear space $(E,\LL)$ is a \underbar{linear space (LS)} if 'at most' in (PS2) is replaced by 'exactly'. We mention in passing that [EMSS] it then holds that $|\LL|\ge |E|$. The equality is e.g. sharp for the Fano plane ($7\ge 7$). 
 A linear space $(E,\Lambda)$ is called\footnote{This kind of 'projective' has nothing to do with projective coverings $(a,b)\approx (c,d)$! Instead of 'projective geometry' one also speaks of 'projective space'.} a \underbar{projective geometry of dimension} $n-1$ if the following takes place.
 
 {\it \begin{enumerate}
 		\item[(PG1)] Every line is long.
 		\item[(PG2)] The Triangle-Axiom (Fig.5.2): If a line intersects two sides of a triangle (not at their intersection), then it also intersects the third line.
 		\item[(PG3)] The  lattice $\La(E,\LL)$ can be spanned by $n$ points but not fewer.
 \end{enumerate}}
 
 \begin{center}
 	\includegraphics[scale=0.7]{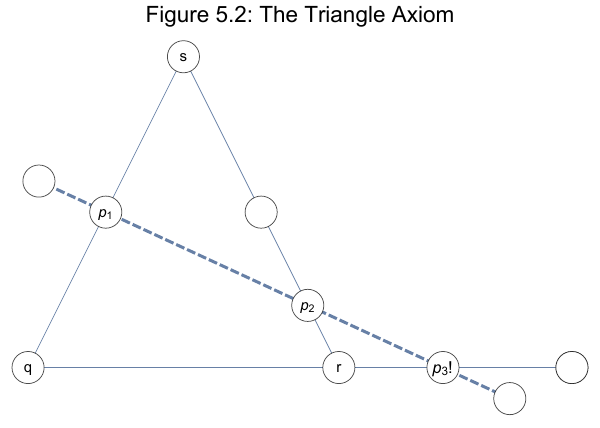} 
 \end{center}
 
 A configuration of four lines as in Fig.5.2 will be called a\underbar{ triangle configuration}.
 In case of PLSes that are projective geometries one speaks of \underbar{subspaces} rather than $\LL$-closed sets.
 Note that (PG3) makes sense for all $n\ge 1$. For $n=1,2$ the projective geometry boils down to
 a single point (so $\Lambda=\es$), respectively to a single line (so $|\Lambda|=1$). For $n=3$ one speaks of a \underbar{projective plane}. For instance the Fano plane is a projective plane.
 
 {\bf 5.6.1} The easiest examples of  projective geometries derive from vector spaces $V$ (over arbitrary skew fields). Namely, let $E_V$ be the set of all 1-dimensional subspaces $\langle p\rangle\ (p\in V)$, and let
 $\LL_V$ be the set of all lines $\ell_S:=\{\langle p\rangle:\ p\in S\}$ where $S$ ranges over 
 all 2-dimensional subspaces. In Ex.5C you are asked to show that $(E_V,\LL_V)$ satisfies (PG1) to (PG3). It is easy to see that $\La(E_V,\Lambda_V)$ is isomorphic to the lattice $\La(V)$ of linear subspaces of $V$. In particular $\La(E_V,\Lambda_V)$ is modular. As is well known, {\it all} projective spaces have modular subspace lattices, but the general proof is more subtle.
 
 In Section 4.5 we considered geometric modular lattices $L$. Let us assume they are c-simple as well. Then if  $\Lambda_L$ is the set of all 'lines' $\ell_x:=At(x)$, where $x$ ranges over the rank two elements, then $(At(L),\Lambda_L)$ turns out to be a projective space as well. It is a classic result that $\La(At(L),\Lambda_L)\simeq L$ (see Ex.5H). More generally, let $(E,\Lambda)$ be {\it any} projective space. Then $L:=\La(E,\Lambda)$ turns out to be a c-simple geometric modular lattice such that $E\simeq At(L)$ and $\Lambda\simeq\Lambda_L$.  
 
 {\bf 5.6.2} In Theorem 7.5 we will drop {\it both} "c-simple" and "geometric" and start investigating arbitrary modular lattices $L$ under a geometric perspective. Specifically, suitable subsets $\ell\s J(L)$ will be called 'lines'. Letting $\Lambda_L$ be the subset of all lines, it was shown by Benson and Conway that $L$ is isomorphic (via $x\mapsto J(x)$) to the closure system of all $\Lambda_L$-closed {\it order ideals} of $(J(L),\le)$. The main perk of Theorem 7.5 is the insight that instead of $\Lambda_L$ it suffices to take suitable (not uniquely determined) {\it subfamilies} $\Lambda\s\Lambda_L$, called 'bases of lines'. Only in 'rare' cases the only choice for $\Lambda$ is $\Lambda=\Lambda_L$, the two most important ones are these. Either $L$ is atomistic (not necessarily c-simple), or $L$ is distributive in which case $\Lambda=\Lambda_L=\es$.
 Indeed, by Birkhoff's Theorem $L\simeq \La(J(L),\le)$ (see Sec.2.5).

 {\bf 5.7} Let $\M(E,cl)$ be a simple matroid and $(J,\Lambda)$ a partial linear space.  Following [W8,Def.4.1] we say that a bijection $\psi:J\to E$ {\it line-preserving (line-pres)} models} $(J,\Lambda)$ if\footnote{It is worth mentioning that some straightforward (yet exponential) algorithm [W8, 4.3.1] can be used to decide whether there is a {\it graphic} matroid that line-pres models $(J,\Lambda)$.} $\psi(\ell)$ is a dependent subset of $E$ whenever $|\ell|\ge 3$; in other words if
  $cl(\{\psi(p),\psi(q)\})=\psi(\ell)$ for all 2-sets $\{p,q\}\s\ell$. If we briefly say that '$\M$ line-pres models $(J,\Lambda)$', then we insinuate that mentioned $\psi$ exists. Upon identifying each $p\in J$ with its 'label' $\psi(p)$ we can view $J$ as a set that simultaneously carries a matroid and a PLS-structure. It then holds that:

\begin{itemize}
	\item[(9)]{\it Each $\M$-closed subset $X\s J$ is $\Lambda$-closed.}
\end{itemize}
 
 Indeed, let $\ell\in\Lambda$ be such that $|\ell\cap X|\ge 2$. Then $\ell=cl(\ell\cap X)\s cl(X)=X$, and so $X$ is $\Lambda$-closed.
 The converse fails; in Fig. 5.3 the set $X=\{a,c,a+c,b+d\}$ is $\Lambda$-closed but not $\M$-closed since $a+b+d\in cl(X)$.
 
 \begin{center}
 	\includegraphics[scale=0.7]{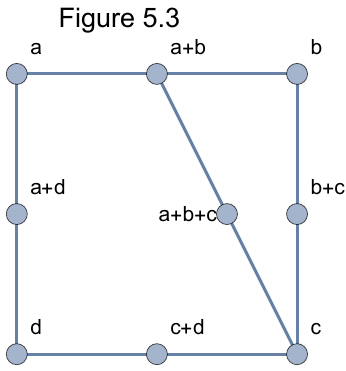} 
 \end{center}
 
 The converse however holds in acyclic PLSes. Specifically (Fig. 5.4(A)) consider a PLS with $J=\{1,2,3\}$ and just a single line, i.e. $\Lambda=\{\{1,2,3\}\}$. Let $k$ be any field with $|k|\ge 3$ and $a,b$ independent $k$-vectors. Defining $\psi(1)=a,\ \psi(2)=a+b,\ \psi(3)=b$ yields a $k$-linear matroid $\M$ on $J$. It is clear that the $\Lambda$-closed subsets of $J$ are exactly the $\M$-closed subsets of $J$. For any further line $\ell$ that we attach to the acyclic PLS so far (such as $\ell=\{3,4,5,6\}$) we invest a new vector not in the span of the previous ones (such as $c$). This vector along with the vector labelling the point of attachment (here $b$) takes care to label the points of $\ell$ (here $\psi(3)=b\ ,\psi(4)=b+c,\ \psi(5)=b+2c,\ \psi(6)=c$). For any acyclic (possibly disconnected) PLS $(J,\Lambda)$ and any large enough\footnote{Akin to $|k|\ge 3$ above we need that $|k|+1$ is $\ge$ the largest cardinality of a line.} field the arising kind of labelling $\psi$ thus line-pres models $(J,\Lambda)$.
 
 \begin{center}
 	\includegraphics[scale=0.7]{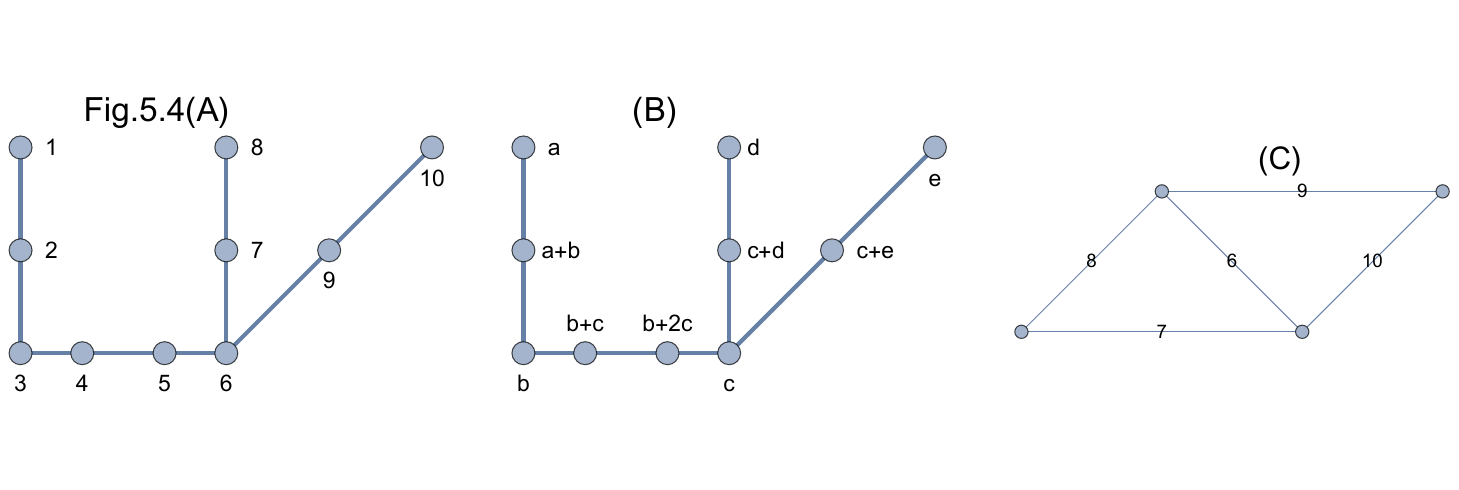} 
 \end{center}
 
 {\bf Theorem 5.2:} {\it Let $(J,\Lambda)$ be an acyclic PLS and let $\M$ be a $k$-linear matroid that  line-pres models $(J,\Lambda)$. Then the $\M$-closed subsets of $J$ coincide with the $\Lambda$-closed subsets of $J$.}
 
 {\it Proof:} Exercise 5F.
 
 For instance in Fig 5.4(B) the set $X=\{a+b,b+2c,c+d,c+e\}$ is $\Lambda$-closed, and so must be $\M$-closed by Theorem 5.2. To spell it out, none of $a,b,b+c,c,d,e$ are in the span of $X$. (This isn't clear for everybody by inspection.)
 
 Following [W3,p.216] we define the \underbar{PLS-rank} of a partial linear space $(J,\Lambda)$ as
 
 $(10)\quad r(J,\Lambda):=|\Lambda|-r^*(J,\Lambda)$
 
 That concept will be crucial in Part C. Among other things we will deal with line-pres $\M$-modelings of $(J,\Lambda)$ that furthermore are \underbar{rank-preserving}  in the sense that the matroid rank $rk(\M)$ equals $r(J,\lambda)+c(J,\Lambda)$. All natural line-pres modelings of acyclic PLSes are rank-preserving (Ex.5G), and certain others as well. For instance the one in Fig.5.3  since $r(J,\Lambda)+c(J,\Lambda)=(5-2)+1=4=rk(\M)$.
 
{\bf 5.7.1}  Let us put the PLSes from Thm. 5.2 in a broader context. Let $\M(E_i,\Ci_i)\ (i=1,2)$ be two matroids given by their circuits which have exactly one point in common, so $E_1\cap E_2=\{p\}$. Define $\Ci:=\Ci_1\cup\Ci_2\cup\Ci^{new}$, where $\Ci^{new}$ consists of all 'novel' circuits $C=C_1\cup C_2\setminus \{p\}$ where
$C_i\in\Ci_i$ and $C_1\cap C_2=\{p\}$. Putting $E:=E_1\cup E_2$ it turns out [Ox,p.240 and Ex.5D] that $\M(E,\Ci)$ is a matroid, the so-called \underbar{parallel connection}
$\M_1*\M_2$ of $\M_1$ and $\M_2$. This construction can be extended to more than two matroids by putting $\M_1*\cdots *\M_n:=(\M_1*\cdots *\M_{n-1})*\M_n$. In fact $*$ turns out (nontrivial) to be associative, and so (say) $\M_1*\cdots*\M_4=\M_1*((\M_2*\M_3)*\M_4)$. Parallel connections tie in nicely [Wh1,Prop.7.6.10] with $k$-linear representations: If $\M_i$ is represented by the set $S_i$ of $k$-vectors, then $\M_1*\cdots *\M_n$ is represented by the set  $S_1\cup\cdots\cup S_n$.

{\bf 5.7.2} We are only concerned with the case where all $\M_i$'s derive from the lines $\ell_i$ of an acyclic PLS. Specifically, we simply declare any three or more points in $\ell_i$ as dependent. Consequently the family $\Ci_i$ of circuits of $\ell_i$ consists of all 3-element subsets of $\ell_i$. For instance, $\ell_2=\{3,4,5,6\}$ in Fig. 5.4(A)
triggers $\M_2(\ell_2,\Ci_2)$ with $\Ci_2=\{456,356,346,345\}$. Exercise 5E asks to calculate $\Ci$ for the parallel connection $\M(E,\Ci)=\M_1*\M_2*\M_3*\M_4$. For instance
$C:=\{1,2,5,7,8\}\in\Ci$. As they must, the five matching vectors in Fig.5.4(B) turn out to be dependent (e.g. $b+2c=((a+b)-a)+2((c+d)-d)$); but any four of them are independent.

{\bf 5.7.3} In [W8] the lines $\ell_i$ of all occuring PLSes have $|\ell_i|\le 3$, this condition being necessary (though not sufficient) for certain associated modular lattices to occur as covering sublattice of a partition lattice. If $|\ell_i|=3$ then we can couple $\ell_i$ to a triangle, more precisely the graphic matroid $\M_i$ induced by the triangle. For instance, if $\ell_3,\ell_4$ are as in Figure 5.4 then $\M_3*\M_4$ is the graphic matroid induced by the graph  in Fig.5.4(C). The term 'parallel connection' (for arbitrary matroids) stems from this  scenario. Specifically, instead of gluing  two graphs by sharing a vertex (serial connection of graphs), here  two graphs are glued by sharing an edge (parallel connection of graphs).

{\bf 5.8 Exercises}

{\bf Exercise 5A.} Prove claim (2).

{\it Solution of Ex.5A.} Let $(E,\LL)$ be a PLS and $S\s\LL$ a connected line-set. Suppose a connected  line-subset $T_1=\{\ell_1,...,\ell_k\}\s S$ has been found (the case $k=1$ is trivial). Let $T_2:=S\setminus T_1$. If $T_2=\es$, we are finished. Otherwise, if the nonvoid pointsets $\bigcup T_1$ and $\bigcup T_2$  were not intersecting, then $(\bigcup S,S)$ would be a disconnected PLS, contradicting the assumption that the line-set $S$ is connected. Hence there is a line $\ell_{k+1}\in T_2$  intersecting $\bigcup T_1$, and so $\{\ell_1,...,\ell_{k+1}\}$ is a larger connected line-set contained in $S$. If the final enumeration is $S=\{\ell_1,...,\ell_n\}$, then $S\setminus\{\ell_n\},\ S\setminus\{\ell_n,\ell_{n-1}\}$, and so forth all remain connected line-sets.

{\bf Exercise 5B.} Let $\M$ be the matroid induced by the projective geometry $(E_V,\LL_V)$ of 5.6.1. Show that
here {\it all} $\LL(\M)$-closed subsets are $\M$-closed. 

{\it Solution of Ex.5B.} Recall from (9) that $\M$-closed are always $\LL(\M)$-closed. As to the converse, if some point $\langle q\rangle$ is in the $\M$-closure of say
$\{\langle p_1\rangle,\langle p_2\rangle,\langle p_3\rangle\}$ then there must be scalars $\lambda_i$ such that $q=\lambda_1 p_1+\lambda_2 p_2+ \lambda_3 p_3$. In order to show that $\langle q\rangle$ is in the $\LL(\M)$-closure of $\{\langle p_1\rangle,\langle p_2\rangle,\langle p_3\rangle\}$ note that   $\langle \lambda_1 p_1+\lambda_2 p_2\rangle$ is in the line determined by $\langle p_1\rangle,\langle p_2\rangle$. Likewise  $\langle (\lambda_1 p_1+\lambda_2 p_2)+\lambda_3 p_3\rangle$ is in the line determined by $\langle \lambda_1 p_1+\lambda_2 p_2\rangle$ and $\langle p_3\rangle$.

{\bf Exercise 5C.}  Let $V$ be a vector space. Show that the coupled linear space $(E_V,\LL_V)$ in 5.6.1 is a projective space.

{\bf Exercise 5D.} Show that $\Ci=\Ci_1\cup\Ci_2\cup\Ci^{new}$ from 5.7.1 satisfies the circuit axioms from 4.2.3.

{\bf Exercise 5E.} Calculate the set family $\Ci$ appearing in the parallel connection $\M(E,\Ci)=\M_1*\M_2*\M_3*\M_4$ in 5.7.2.

{\it Solution of Ex.5E.} The matroid $\M_1(\ell_1,\Ci_1)$ coupled to line $\ell_1$ has $\Ci_1=\{\ell_1\}=\{\{1,2,3\}\}\ '='\ \{123\}$ (for brevity). We likewise have
$\M_i(\ell_i,\Ci_i)$ with $\Ci_2=\{456,356,346,345\},\ \Ci_3=\{678\},\ \Ci_4=\{6910\}\ (=\{6,9,10\})$.  If $\M'(E',\Ci')$ is the parallel connection of $\M_1$ and $\M_2$,
$\M''(E'',\Ci'')$ is the parallel connection of $\M'$ and $\M_3$, and $\M(E,\Ci)$ is the parallel connection of $\M''$ and $\M_4$, then one checks that
\begin{itemize}
	\item $\Ci'=\Ci_1\cup\Ci_2\cup\{1256,1246,1245\}$
	\item $\Ci''=\Ci'\cup\Ci_3\cup\{4578,3578,3478,12578,12478\}$
	\item $\Ci=\Ci''\cup\Ci_4\cup\{45910,3510,3410,78910\}=\Ci_1\cup\cdots\cup\Ci_4\cup\{1256,\ldots,78910\}$
	\end{itemize}

{\bf Exercise 5F.} Prove Theorem 5.2.

{\bf Exercise 5G.} Show that the natural line-pres modelling of every acyclic PLS (by a $k$-linear or graphic  matroid) is rank-preserving.

{\bf Exercise 5H.} Let $L$ be an atomistic modular lattice (not necessarily c-simple). As for any lattice, $L$ is isomorphic to the lattice (closure system) $\{J(x):\ x\in L\}$; see Sec.2.4 and Ex.2C. In order to prove $(At(L),\Lambda_L)\simeq L$ (see 5.6.1), it thus suffices to show that the $\Lambda_L$-closed subsets of $At(L)$ are exactly the sets of type $J(x)$. Therefore show:
\begin{enumerate}
	\item[(a)] The set $J(x)\s At(L)$ is $\Lambda_L$-closed for all $x\in L$ (easy).
	\item[(b)] If $X\s At(L)$ is $\Lambda_L$-closed, then $X=J(x)$ for $x:=\sum X$ (more subtle). 
\end{enumerate}

\bigskip

 \centerline{\bf 6.   Enumeration of modular lattices }
 
 Let $g(n)$ be the number of (non-isomorphic) distributive lattices of cardinality $n$. According to [EHK] it holds that $g(30)=8'186'962$, the largest known value being $g(49)=908'414'736'485$.
 For counting modular lattices completely different methods are in order.
 Kohonen [K] has determined the number $f(n)$ of modular, non-distributive $n$-element lattices for all $n\le 30$, thereby shattering the old record of $n=23$. For instance $f(30)=3'485'707'007 >>g(30)$.
 
  A lattice is \underbar{vertically decomposable } if there is $x\in L\setminus\{0,1\}$ such that for all $y\in L$ one has $y\le x$ or $y\ge x$. Then $L$ is the 'vertical sum' of $L_1:=[0,x]$ and $L_2:=[x,1]$. Conversely, given any lattices $L_1',L_2'$ it is easy (Ex.6A) to construct a lattice $L$ wich is the vertical sum of sublattices $L_1\simeq L_1'$ and $L_2\simeq L_2'$. In both cases we write $L=L_1\oplus L_2$. For instance, Figure 6.1 (A) shows $M_3\oplus M_4$. Vertically decomposable lattices reduce in a trivial manner to their constituents and whence are often neglected in the enumeration of lattices. More relevant than $M_3\oplus M_4$ is $M_{3,4}$ (Fig. 6.1(B)). Generally
 $M_{s,t}$ is obtained by gluing\footnote{Such gluings will be discussed in Section 13. They are not related to (matroidal) parallel connections.}  an upper covering of $M_s$ and a lower covering of $M_t$. 
 In  Section 6 it is convenient to allow $n=2$ in $M_n$, so $M_2$ is the distributive length 2 lattice with 2 atoms; see $M_{3,2}$ in Figure (C) and $M_{2,2}$ in (D).
 
 \begin{center}
 	\includegraphics[scale=0.7]{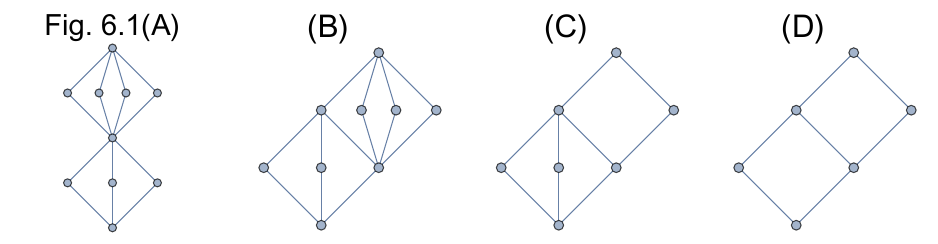} 
 \end{center}

 \bigskip

 Jukka Kohonen gratefully provided the diagrams of all 56 vertically indecomposable modular lattices of cardinality at most 10, the distributive ones among them (say $M8.5^*$) have the superscript $\ast$:

 \begin{center} \includegraphics[scale=0.35]{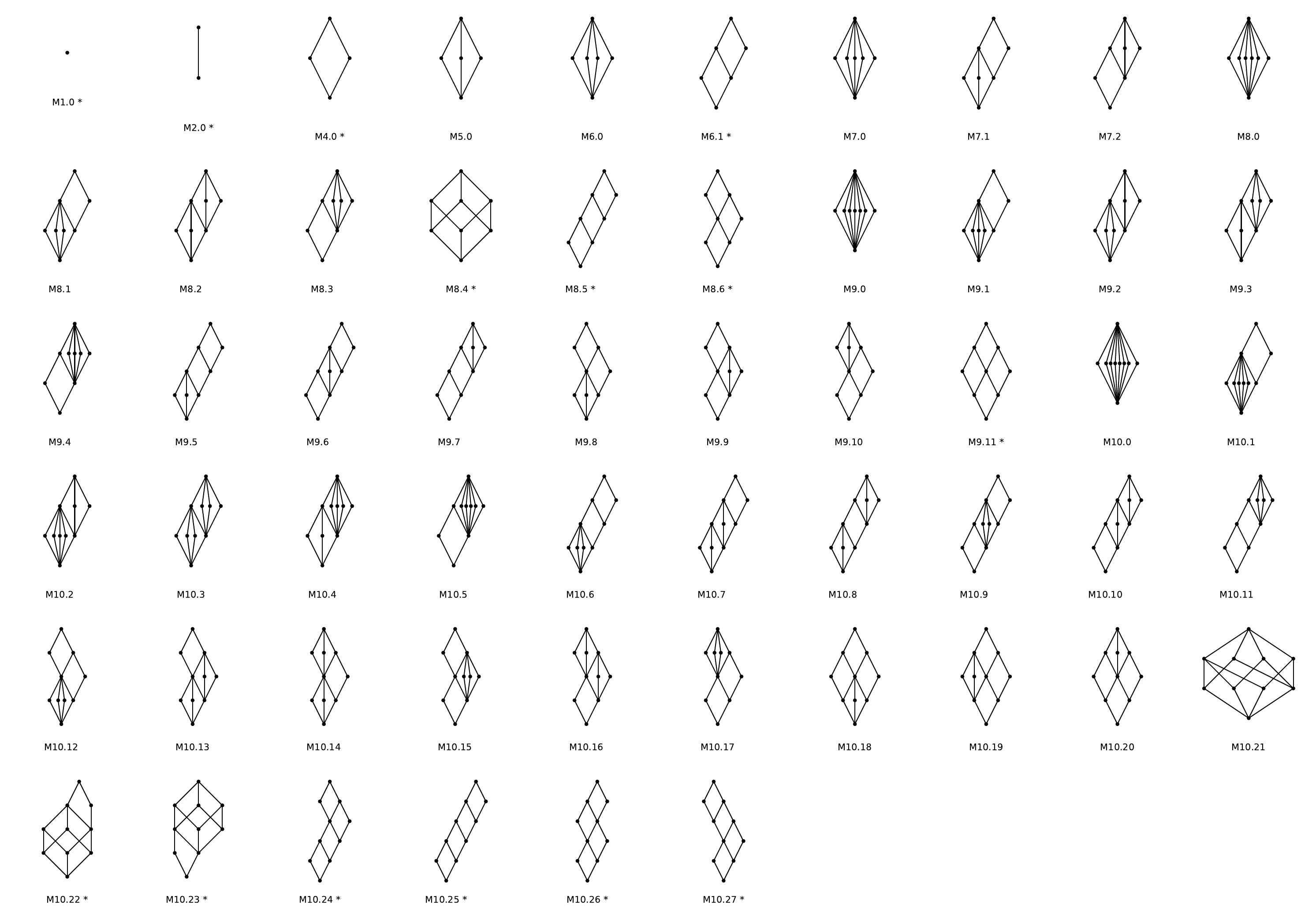} 
\end{center} 

And here are the 54 vertically indecomposable modular lattices with exactly 11 elements:

\begin{center} \includegraphics[scale=0.35]{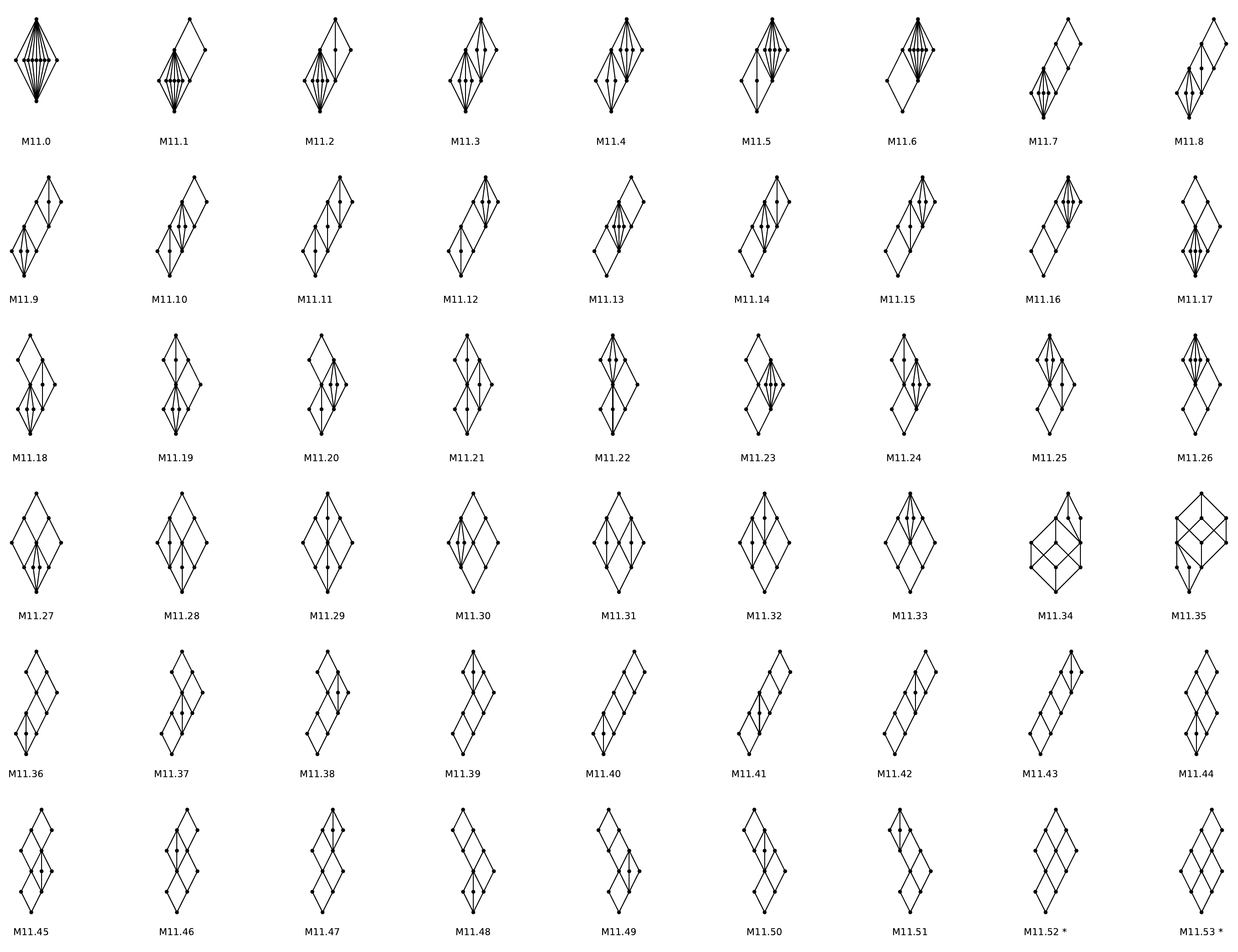} 
\end{center}

 Diagrams of all modular lattices up to cardinality 16 can be requested per email; either write to mwild@sun.ac.za or to jukka.kohonen@aalto.fi.
 
 {\bf 6.1} Here comes a fairly easy-to-implement algorithm to check the congruence-simplicity of any modular lattice $L$ whose cover relations $a\prec b$ are known. (The reasons for WHY this works will be provided in Part B.) First, identify $J(L)$, i.e. all elements with exactly one lower cover. Second,  partition $L$ into classes of elements of equal rank. Third, find all \underbar{$M_t$-elements} $x$ as follows. For all $x\in L\setminus J(L)$ find all its $t$ lower covers $x_1,x_2,..$. If $t= 2$, go to the next $x$. Otherwise check whether in the class of all elements of rank $\delta(x)-2$ there is some (necessarily unique) $x_0$ which is below all $x_i$'s. If no, go to the next $x$. If yes, add the $M_t$-element $x$ to the set of $M_t$-elements and check the next $x$.
 Fourth, for each $M_t$-element $x$ find any 'line' $\ell_x=\{p_1,...,p_t\}$ by picking (at random) $p_i$ from $J(x_0,x_i)$.  Fifth, check whether the PLS $(J(L),\{\ell_x:\ x\ is\ M_t$-$element\})$ is connected. If yes, $L$ is c-simple, otherwise not.  One can speed up the algorithm by doing the following between Step 1 and Step 2. Check whether any $p\in J(L)$ is join-prime (using Theorem 2.1). If yes, $L$ is not c-simple. If no, continue with Steps 2 to 5.

{\bf 6.2 Exercises}

{\bf Exercise 6A.} Given any lattices $L_1',L_2'$ carefully define a lattice $L$ which is the vertical sum of sublattices $L_1\simeq L_1'$ and $L_2\simeq L_2'$.


{\bf Exercise 6B.} It turns out (by brute-force inspection of the Kohonen diagrams) that there are exactly 11 congruence-simple modular lattices $L$ with $9\le |L|\le 12$. They can be obtained by adding doubly irreducible elements to the three lattices below.
Please do so.

\begin{center}
	\includegraphics[scale=0.7]{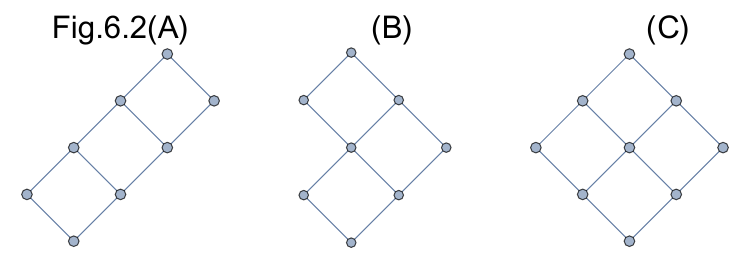} 
\end{center}

\bigskip
\centerline{\bf References (to be updated)}
\medskip

\begin{itemize}

\item[[B]] G. Birkhoff, Lattice Theory, AMS Coll. Publ., vol.25, third
edition 1967.
\item[[BB]] L.M. Batten, A. Beutelsbacher, The theory of finite linear spaces, Cambridge University Press 2009.
\item[[BC]] D. Benson, J. Conway, Diagrams for modular lattices, J. Pure Appl. Algebra 37 (1985) 111-116.
\item[[CD]] J. Crawley and R.P. Dilworth, Algebraic Theory of Lattices,
Prentice Hall, Englewood Cliffs 1973.

\item[[DF]] A. Day and R. Freese, The role of glueing constructions in 
modular lattice theory, p.251-260 in: The Dilworth Theorems, Birkhauser 1990.
\item[[DP]] BA Davey, HA Priestley, Introduction to lattices and order, Cambridge University Press 2002.
\item[[EHK]]  Marcel Ern\'e, Jobst Heitzig, and J\"urgen Reinhold, On the number of distributive lattices,
Electron. J. Combin. 9 (2002).
\item[[EMSS]]  Erdős, P; Mullin, R. C.; Sós, V. T.; Stinson, D. R.;
Finite linear spaces and projective planes.
Discrete Math. 47 (1983), no. 1, 49–62.
\item[[FH]] U. Faigle and C. Herrmann, Projective geometry on partially
ordered sets, Trans. Amer. Math. Soc. 266 (1981) 319-332.
\item[[FRS]] U. Faigle, G. Richter and M. Stern, Geometric exchange
properties in lattices of finite length, Algebra Universalis 19 (1984)
 355-365.
\item[[G]] G. Gr\"atzer, General Lattice Theory: Foundation, Birkh\"auser 2011.
\item[[H]] C. Herrmann, S-verklebte Summen von Verbaenden, Math. Zeitschrift 130 (1973) 225-274.
\item[[HPR]] C. Herrmann, D. Pickering and M. Roddy, Geometric description
of modular lattices, to appear in Algebra Universalis.
\item[[HW]] C. Herrmann and M. Wild, Acyclic modular lattices and their
representations, J. Algebra 136 (1991) 17-36.
\item[[Hu]] A. Huhn, Schwach distributive Verb\"ande I, Acta Sci. Math. Szeged
33 (1972) 297-305.
\item[[JN]] B. J\'onsson and J.B. Nation, Representations of $2$-distributive
modular lattices of finite length, Acta Sci. Math. 51 (1987) 123-128.
\item[[K]] J. Kohonen, Generating modular lattices of up to 30 elements. Order 36 (2019) 423–435.
\item[[KNT]] McKenzie, Ralph ; McNulty, George ; Taylor, Walter F;
Algebras, lattices, varieties. Vol. I.
 Wadsworth and Brooks/Cole Advanced Books and Software, Monterey, CA, 1987. xvi+361 pp. 
\item[[LMR]] 
Lux, Klaus; Müller, Jürgen; Ringe, Michael;
Peakword condensation and submodule lattices: an application of the MEAT-AXE. 
J. Symbolic Comput. 17 (1994), no. 6, 529–544.
\item[[P]] W. Poguntke,  Zerlegung von S-Verbänden, Math. Z. 142 (1975), 47–65.
\item[[S]] M. Stern, Semimodular lattices, Teubner Texte zur Mathematik 125
(1991), Teubner Verlagsgesellschaft Stuttgart, Leipzig, Berlin.
\item[[Wh1]] N. White (ed.), Theory of Matroids, Encyclopedia Math. and Appl.
26, Cambridge University press 1986.
\item[[Wh2]] N. White (ed.), Combinatorial geometries,
 Encyclopedia Math. and Appl. 29, Cambridge University press 1987.

\item[[W1]] M. Wild, Dreieckverb\"ande: Lineare und quadratische 
Darstellungstheorie, Thesis University of Zurich 1987.
  
\item[[W2]] M. Wild, Modular lattices of finite length, first version (1992) of the present paper, ResearchGate.
\item[[W3]] M. Wild, Cover preserving embedding of modular lattices into
partition lattices, Discrete Mathematics 112 (1993) 207-244.
\item[[W4]] M. Wild, The minimal number of join irreducibles of a finite modular lattice. Algebra Universalis 35 (1996),  113–123.

\item[[W5]] M. Wild, Optimal implicational bases for finite modular lattices. Quaest. Math. 23 (2000), no. 2, 153–161.
\item[[W6]] M. Wild, Output-polynomial enumeration of all fixed-cardinality ideals of a poset, respectively all fixed-cardinality subtrees of a tree, Order 31 (2014) 121–135.
\item[[W7]] M. Wild, The joy of implications, aka pure Horn formulas: mainly a survey. Theoret. Comput. Sci. 658 (2017), part B, 264–292.
\item[[W8]] M. Wild, Tight embedding of modular lattices into partition lattices: progress and program, Algebra Universalis 79 (2018) 1-49.

\item[[Wi1]] R. Wille, Über modulare Verbände, die von einer endlichen halbgeordneten Menge frei erzeugt werden. 
Math. Z. 131 (1973), 241–249.
\item[[Wi2]] Subdirekte Produkte vollständiger Verbände. (German) J. Reine Angew. Math. 283(284) (1976), 53–70.
\end{itemize}

\end{document}